\documentclass[11pt]{article}
\usepackage{amssymb,epsfig,amsmath}
\usepackage[dvips]{color}
\usepackage[T1]{fontenc}
\usepackage[francais]{babel}
\usepackage[utf8]{inputenc}
\usepackage{comment}
\usepackage{hyperref}
\usepackage{graphicx}
\usepackage{multicol}
\usepackage{caption}
\newtheorem{defi}{D\'efinition}[section]
\newtheorem{prop}[defi]{Proposition}
\newtheorem{theo}[defi]{Th\'eor\`eme}
\newtheorem{conj}[defi]{Conjecture}
\newtheorem{lemm}[defi]{Lemme}
\newtheorem{coro}[defi]{Corollaire}
\newtheorem{rema}[defi]{Remarque}
\newtheorem{exem}[defi]{Exemple}
\newtheorem{exems}[defi]{Exemples}

\newcommand{\bdefi}{\begin{defi}}
\newcommand{\edefi}{\end{defi}}
\newcommand{\bprop}{\begin{prop}}
\newcommand{\eprop}{\end{prop}}
\newcommand{\btheo}{\begin{theo}}
\newcommand{\etheo}{\end{theo}}
\newcommand{\blemm}{\begin{lemm}}
\newcommand{\brema}{\begin{rema}}
\newcommand{\erema}{\end{rema}}
\newcommand{\bexer}{\begin{exem}}
\newcommand{\eexer}{\end{exem}}
\newcommand{\bexems}{\begin{exems}}
\newcommand{\eexems}{\end{exems}}
\newcommand{\bconj}{\begin{conj}}
\newcommand{\econj}{\end{conj}}
\newcommand{\elemm}{\end{lemm}}
\newcommand{\bcoro}{\begin{coro}}
\newcommand{\ecoro}{\end{coro}}
\newcommand{\dem}{\noindent{\bf Démonstration. }}
\newcommand{\rem}{\noindent{\bf Remarque. }}

\usepackage{mathrsfs}
\renewcommand\mathcal{\mathscr}

\newcommand{\T}{{\cal T}}

\newcommand{\M}{{\cal M}}

\newcommand{\G}{{\cal G}}

\newcommand{\HHHH}{{\cal H}}
\newcommand{\LLLL}{{\cal L}}
\newcommand{\PPPP}{{\cal P}}
\newcommand{\C}{{\cal C}}

\newcommand{\maths}[1]{{\mathbb #1}}  

\newcommand{\RR}{\maths{R}}
\newcommand{\NN}{\maths{N}}

\newcommand{\ZZ}{\maths{Z}}

\newcommand{\cqfd}{\hfill$\Box$}

\newcommand{\CAT}{\operatorname{CAT}}

\newcommand{\dddp}{\partial_{\infty}^2\widetilde{\Sigma}}
\newcommand{\ddp}{\partial_{\infty}\widetilde{\Sigma}}

\newcommand{\Image}{\operatorname{Image}}

\newcommand{\Supp}{\operatorname{Supp}}

\newcommand{\revet}{(\widetilde{\Sigma},[\widetilde{q}])}

\newcommand{\srfce}{({\Sigma},{[q]})}

\newcommand{\revetm}{(\widetilde{\Sigma},\widetilde{m})}
\newcommand{\srfcem}{({\Sigma},m)}
\newcommand{\Id}{\operatorname{Id}}

\newcommand{\Lq}{\Lambda_{[q]}}

\newcommand{\Lqr}{\widetilde{\Lambda}_{[\widetilde{q}]}}

\newcommand{\Lmr}{\widetilde{\Lambda}_{\widetilde{m}}}
\newcommand{\Gqr}{\G_{[\widetilde{q}]}}
\newcommand{\Gmr}{\G_{\widetilde{m}}}
\newcommand{\Gq}{\G_{[q]}}

\newcommand{\Ltilde}{\widetilde{\Lambda}}
\newcommand{\mutilde}{\widetilde{\mu}}
\newcommand{\grperevet}{\Gamma_{\widetilde{\Sigma}}}

\newcommand{\BP}{\operatorname{BP}}
\title{Laminations géodésiques plates mesurées.}
\author{Thomas~Morzadec}

\begin{document}
\maketitle

Département de mathématique, UMR 8628 CNRS, Université Paris-Sud, Bât. 430, F-91405 Orsay Cedex, France  
Bureau : 16.  {\it thomas.morzadec@math.u-psud.fr}

\medskip

\textbf{Abstract:} In \cite{Morzy1}, we have introduced a notion of flat laminations on surfaces endowed with a flat structure, similar to
geodesic laminations on hyperbolic surfaces.
Here is a sequel to this article that aims at defining transversal measures on flat laminations similar to
transversal measures on hyperbolic laminations, taking into account that two different leaves of a flat lamination may no longer disjoint.  
Then, we define a topology on the set of  measured flat laminations and show that its projective space is compact. Finally, we define the dual tree to a 
measured flat lamination.   
\footnote{ Mots clés: Lamination géodésique mesurée, surface munie d'une structure plate à holonomie $\{\pm\Id\}$, différentielle quadratique holomorphe, 
feuilletage mesurée, surface hyperbolique, arbre dual.
Code AMS 30,37,53,57,58.}

\medskip

\textbf{Résumé:} Dans \cite{Morzy1}, nous avons introduit une notion de laminations sur les surfaces munies d'une structure plate,  analogues aux laminations
géodésiques sur les surfaces hyperboliques.
Cet article est une suite de \cite{Morzy1} dont le but est de définir les mesures transverses sur les laminations plates analogues aux mesures transverses
sur les laminations hyperboliques, en prenant garde que deux feuilles distinctes ne sont plus nécessairement disjointes.  
Ensuite, nous définissons une topologie sur l'ensemble des laminations plates mesurées et nous montrons que son projectifié est compact.
Enfin, nous définissons l'arbre dual à une lamination plate mesurée.

\section{Introduction.}

Dans \cite{Morzy1}, nous avons proposé et étudié une définition des laminations géodésiques sur les espaces métriques localement $\CAT(0)$ enrubannés, 
analogues aux laminations géodésiques sur les surfaces hyperboliques, en particulier sur les
surfaces munies de structures plates à 
holonomie $\{\pm\Id\}$, que nous avons appelées {\it laminations plates}. Dans cet article, nous proposons une définition de {\it lamination plate mesurée} 
  et nous étudions le lien entre laminations plates mesurées et laminations hyperboliques mesurées (voir par exemple \cite{Bonahon97}).
             Nous reprendrons les notations de \cite{Morzy1}, notamment nous noterons $\srfce$ une surface munie d'une structure plate
(à holonomie $\{\pm\Id\}$) et 
$\Lambda$ une lamination plate sur $\srfce$.

La principale difficulté par rapport aux mesures transverses sur les laminations hyperboliques est que les feuilles ne sont pas disjointes en général,
donc on ne va pas définir des mesures sur les images des arcs transverses à la lamination, mais sur l'ensemble des géodésiques locales qui les intersectent 
transversalement,   
     et il faut redéfinir la notion d'invariance par holonomie de ces familles de mesures.

Dans la première partie, nous définissons les laminations plates mesurées et munissons leur ensemble d'une topologie.
Dans la deuxième, nous définissons
l'image réciproque d'une lamination plate mesurée dans un revêtement. Dans la troisième, nous définissons un homéomorphisme entre 
l'espace des laminations plates mesurées sur une surface (compacte) munie d'une structure plate 
et l'espace des mesures de Radon sur l'ensembles des géodésiques (définies à changements d'origine près) d'un de ses revêtements universels
(localement isométriques),
qui sont invariantes par l'action du groupe de revêtement et dont le support est une lamination plate.    
%
        Dans la quatrième,          
nous établissons une application continue, surjective et propre des laminations plates mesurées dans les laminations hyperboliques mesurées, pour la  métrique 
hyperbolique complète dans la classe conforme. Dans la cinquième, nous définissons le nombre d'intersection entre une lamination plate mesurée et une classe
d'homotopie libre de courbes fermées. 
Dans la dernière partie, nous définissons l'arbre dual à une lamination plate mesurée, et l'action du groupe de revêtement universel sur cet arbre.

Je tiens à remercier chaleureusement Frédéric Paulin pour ses nombreux conseils et ses relectures attentives.

\section{Rappels et conventions.}\label{rappelsconventions}

Dans tout ce texte, nous utiliserons les définitions et les notations de \cite{Morzy1} concernant les surfaces
munies de structures plates (à holonomie $\{\pm\Id\}$) et les laminations plates (voir notamment les définitions \cite[Déf.~2.2]{Morzy1} et
\cite[Déf.~2.8]{Morzy1}).
Sauf mention du contraire, nous noterons
$\srfce$ une surface connexe, orientable, à bord (éventuellement vide) munie d'une structure plate et $p:\revet\to\srfce$ un revêtement universel localement isométrique.
Si $\Sigma$ est compacte, on supposera toujours que $\chi(\Sigma)<0$. 
     On munira le bord à l'infini $\partial_\infty\widetilde{\Sigma}$ de l'ordre cyclique défini par le choix d'une orientation de
$\Sigma$ (la surface $\Sigma$ sera toujours compacte ou un revêtement d'une surface compacte, voir la remarque \cite[Rem.~2.9]{Morzy1}), et on notera 
$\dddp=\ddp\times\ddp-\{(x,x),x\in\ddp\}$.
On rappelle que la {\it
topologie des géodésiques} est la topologie compacte-ouverte sur l'ensemble $\G_d$ des géodésiques locales paramétrées pour la distance $d$ (noté $\G_{[q]}$
pour une structure plate $[q]$ et $\G_m$ pour une métrique hyperbolique $m$) ou le quotient de la topologie compacte-ouverte
 par l'action par translations à la source de $\RR$, sur l'ensemble $[\G_d]$ des géodésiques locales définies à changements d'origine près.
 Ce seront les seules topologies considérées sur les 
 espaces de géodésiques locales. On notera $g\mapsto[g]$ l'application de passage au quotient de $\G_d$ dans $[\G_d]$, et si $F\subset\G_d$, on notera
 $[F]$ son image dans $[\G_d]$. On appelle {\it arc} une application $\C^1$ par morceaux $\alpha:[0,1]\to\Sigma$ qui est un homéomorphisme sur
 son image.

\medskip

        Soit $\Lambda$ une 
lamination plate de $\srfce$.  
Un arc $\alpha$ est {\it transverse} à une feuille ou segment de feuille $\ell$ de $\Lambda$ si

\noindent
$\bullet$~$\alpha$ est transverse à $\ell$ en dehors des singularités de $[q]$ et des points singuliers de $\alpha$;
\medskip
\noindent
\begin{minipage}{12 cm}
\noindent
$\bullet$~ pour toute singularité $x$ de $[q]$ ou point singulier de $\alpha$ appartenant à $\Image(\ell)\cap\alpha(]0,1[)$, 
il existe un voisinage
$D$ de $x$ qui est un disque topologique, et
un segment $S$ de $\ell$ tels que $D-\Image(S)\cap D$ a deux composantes connexes et les composantes connexes de $D\cap(\alpha([0,1])-\{x\})$ sont contenues
dans des composantes distinctes de $D-\Image(S)\cap D$;
\end{minipage}
\begin{minipage}{2.9 cm}
\begin{center}
\begin{picture}(0,0)%
\includegraphics{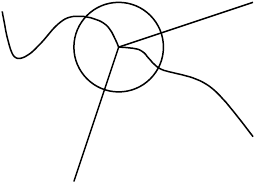}%
\end{picture}%
\setlength{\unitlength}{3771sp}%
\begingroup\makeatletter\ifx\SetFigFont\undefined%
\gdef\SetFigFont#1#2#3#4#5{%
  \reset@font\fontsize{#1}{#2pt}%
  \fontfamily{#3}\fontseries{#4}\fontshape{#5}%
  \selectfont}%
\fi\endgroup%
\begin{picture}(1284,924)(1204,-1873)
\put(1801,-1096){\makebox(0,0)[lb]{\smash{{\SetFigFont{11}{13.2}{\familydefault}{\mddefault}{\updefault}{\color[rgb]{0,0,0}$x$}%
}}}}
\put(1486,-1501){\makebox(0,0)[lb]{\smash{{\SetFigFont{11}{13.2}{\familydefault}{\mddefault}{\updefault}{\color[rgb]{0,0,0}$U$}%
}}}}
\put(2116,-1231){\makebox(0,0)[lb]{\smash{{\SetFigFont{11}{13.2}{\familydefault}{\mddefault}{\updefault}{\color[rgb]{0,0,0}$\ell$}%
}}}}
\put(2071,-1501){\makebox(0,0)[lb]{\smash{{\SetFigFont{11}{13.2}{\familydefault}{\mddefault}{\updefault}{\color[rgb]{0,0,0}$\alpha$}%
}}}}
\end{picture}%

\end{center}
\end{minipage}

\noindent
$\bullet$~ $\alpha$ n'est tangent à $\Image(\ell)$ ni en $0$ ni en $1$.

\medskip

On dit que $\alpha$ est {\it transverse à un ensemble $F$ de feuilles} ou de segments de feuilles de $\Lambda$ si $\alpha$ est transverse à chacun des éléments
de $F$, et que $F$ est {\it transverse à $\alpha$} si $\alpha$ est transverse à $F$. Notamment,
un arc est transverse à $\Lambda$ s'il est transverse à chacune des feuilles de $\Lambda$.  

   Si $\alpha:[0,1]\to\Sigma$ est 
un arc de $\Sigma$, on note $G(\alpha)$ l'ensemble des éléments de $\Gq$ transverses à $\alpha$ dont l'origine appartient
à $\alpha([0,1])$. Soit $F_1\subseteq\Gq$ tel que $[F_1]\subseteq\Lambda$ et soient $\alpha_1$ et $\alpha_2$ deux arcs disjoints, transverses à $F_1$,
tels que $F_1\subseteq G(\alpha_1)$ et tout élément de $F_1$ intersecte $\alpha_2([0,1])$ en un temps positif. Pour tout $g_1\in F_1$, on note $t_{g_1}=
\min\{t>0\;:\;g_1(t)\in \alpha_2([0,1])\}$, et $F_2$ le sous-ensemble des éléments $g_2\in G(\alpha_2)$ tels qu'il existe $g_1\in F_1$ telle que 
$g_2(t)=g_1(t+t_{g_1})$ pour tout $t\in\RR$. Une {\it holonomie} $h:F_1\to F_2$ de $\Lambda$
est un homéomorphisme entre $F_1$ et $F_2$ défini par $h(g_1)=g_2:t\mapsto g_1(t+t_{g_1})$ tel qu'il existe une homotopie
$H:[0,1]\times[0,1]\to\Sigma$ entre $\alpha_1$ et $\alpha_2$ telle que :

\noindent
$\bullet$~ pour tout $t\in\,[0,1]$, l'application $s\mapsto H(s,t)$ est un arc transverse à chacun des segments de feuilles $g_{1|[0,\,t_{g_1}]}$, pour $g_1\in
F_1$;

\noindent
$\bullet$~ pour tout $\ell\in F_1$, il existe $s_\ell\in[0,1]$ tel que $t\mapsto H(s_\ell,t)$ soit un segment de $\ell$ (à reparamétrage près);

\noindent
$\bullet$~ les intersections $H([0,1]\times]0,1[)\cap \alpha_i([0,1])$ pour $i=1,2$ sont vides.

\bdefi\label{defmesuretransverse} 
Une mesure transverse à $\Lambda$ est la donnée pour tout arc $\alpha$ transverse à $\Lambda$ d'une mesure de Radon $\mu_\alpha$ sur $G(\alpha)$ telle que:

\noindent
$(1)$~ Le support de $\mu_\alpha$ est l'ensemble $\{\ell\in G(\alpha)\;:\;[\ell]\in\Lambda\}$;

\noindent
$(2)$~ si $h:F_1\to F_2$ est une holonomie de $\Lambda$, avec $F_1\subset G(\alpha_1)$ et $F_2\subset G(\alpha_2)$, où $F_1$ est un borélien de $G(\alpha_1)$
et $\alpha_1$ et $\alpha_2$
sont deux arcs transverses à $F_1$, alors $h_*(\mu_{\alpha_1|F_1})=\mu_{\alpha_2|F_2}$;

\noindent
$(3)$~ $\mu_\alpha$ est $\iota$-invariante, où $\iota(\ell)=\ell^-:t\mapsto\ell(-t)$;

\noindent
$(4)$~ si $\alpha'([0,1])\subseteq\alpha([0,1])$, alors $\mu_{\alpha|G(\alpha')}=\mu_{\alpha'}$. 
\edefi

On notera $(\Lambda,\mu)$ une lamination plate munie d'une mesure transverse, que l'on appellera {\it lamination plate mesurée}, et $\M\LLLL_p(\Sigma)$
l'espace des 
laminations plates mesurées sur $\Sigma$.
%
   On munit $\M\LLLL_p(\Sigma)$ de la topologie telle qu'une suite $(\Lambda_n,\mu_n)_{n\in\NN}$ converge vers $(\Lambda,\mu)$ si et seulement si pour tout
arc $\alpha$
transverse à $\Lambda$, $\alpha$ est transverse à $\Lambda_n$ pour $n$ assez grand et $\mu_{n,\alpha}\overset{*}{\rightharpoonup}\mu_\alpha$ dans l'espace des mesures de 
Radon sur $G(\alpha)$.

On dit qu'une feuille $\ell$ de $\Lambda$ est {\it positivement récurrente} si on peut fixer son origine en dehors d'une singularité, 
de sorte que pour tout $T,\varepsilon>0$,
il existe $t>T$ tel que $d(\ell(t),\ell(0))<\varepsilon$. Par exemple, si
$\Sigma$ est compacte, les feuilles qui aboutissent (au sens de \cite[Th.~6.1]{Morzy1})  dans des composantes cylindriques ou minimales sont positivement
récurrentes.   

\blemm\label{quedescomposnatesminimales}
Si $\Lambda$ est munie d'une mesure transverse $\mu$, alors $\Lambda$ n'a pas de feuille isolée 
    qui soit positivement récurrente sans être
périodique. 
\elemm

\medskip

\dem Supposons qu'il existe une feuille isolée $\ell$ de $\Lambda$ qui soit positivement récurrente, dont on fixe l'origine telle que ci-dessus.
      Alors il existe un segment géodésique $\alpha$ de diamètre assez petit, ne rencontrant pas de singularité,
 transverse à $\Lambda$, qui intersecte $\ell$ en son origine, et une suite strictement croissante 
 $(t_n)_{n\in\NN}$ des instants positifs successifs où $\ell$ intersecte l'image de $\alpha$.   
 Pour tout $n\in\NN$, soit $\ell_n:t\mapsto\ell(t+t_n)$. 
 Quitte à restreindre $\alpha$ et à changer sa direction, on définit une homotopie $H:[0,1]\times[0,1]\to\Sigma$ entre $\alpha$ et lui même en le translatant
 parallèlement le long de chacune des liaisons de singularités parcourues par $\ell_{|[t_n,t_{n+1}]}$, et en prolongeant
 aux singularités pour que $H$ soit continue, de sorte que tous les arcs $s\mapsto H(s,t)$, avec $t\in[0,1]$, soient transverses à $\ell_{|[t_n,t_{n+1}]}$ 
 (cette construction n'est pas unique). Puisque $[t_n,t_{n+1}]$ est compact, quitte à restreindre $\alpha$, 
 il existe un arc $\beta$ tel qu'il existe un unique  $z\,\in]0,1[$ tel que $\beta:=s\mapsto H(s,z)$ et $\beta([0,1])\cap H([0,1]\times
 ([0,z[\cup]z,1]))=\emptyset$. Soit
 $t_\beta$  l'unique réel de $]t_n,t_{n+1}[$ tel que $\ell(t_\beta)\in\beta([0,1])$. Alors il existe une holonomie entre $\{\ell_n\}$ et $\{\ell_{n,\beta}\}$,
 où $\ell_{n,\beta}:t\mapsto\ell(t+t_\beta)$ appartient à $G(\beta)$. De même, il existe une holonomie entre
 $\{\ell_{n,\beta}\}$ et $\{\ell_{n+1}\}$. D'après le point $(2)$ de la définition \ref{defmesuretransverse}, on en déduit que 
 $\mu_\alpha(\{\ell_n\})=\mu_\alpha(\{\ell_{n+1}\})$, et par récurrence, $\mu_\alpha(\{\ell_n\})=\mu_\alpha(\{\ell_0\})$. Comme $\ell$ est isolée,   
      chacune des feuilles $\ell_n$ est isolée dans $G(\alpha)$, et puisque $\ell_0$ appartient au support de $\mu_\alpha$, on a $\mu_\alpha(\ell_0)>0$. Or, 
      l'ensemble $\{\ell_n\}_{n\in\NN}$ est contenu
 dans $G(\alpha)$, qui est relativement compact d'après le théorème d'Ascoli. Mais, puisque $\ell$ n'est pas périodique, cet ensemble 
 est infini et $\mu_\alpha(\overline{\{\ell_n\}_{n\geqslant 0}})
 \geqslant\sum_{n\geqslant 0}\mu_\alpha(\ell_0)=+\infty$. Mais alors $\mu_\alpha$ ne serait pas localement finie.\cqfd
 
 \section{Relevé d'une lamination plate mesurée.}\label{relevez}
 
 On note toujours $\srfce$ une surface munie d'une structure plate comme dans la partie \ref{rappelsconventions} et $(\Lambda,\mu)$ une lamination plate mesurée
 sur $\Sigma$. 
 Soit $p':(\Sigma',[q'])\to\srfce$ un revêtement localement isométrique de $\srfce$ de groupe de revêtement
 $\Gamma_{\Sigma'}$ et $\Lambda'$ l'image 
 réciproque de $\Lambda$ dans $\Sigma'$ (voir la remarque précédant le lemme \cite[Lem.~2.3]{Morzy1}). Puisque $p'$ est un difféomorphisme local, 
 si $\alpha$ est un arc de $\Sigma$ transverse à $\Lambda$, et si $\alpha'$ est un relevé de $\alpha$ dans $\Sigma'$, alors $\alpha'$ est transverse à $\Lambda'$,
 et $p'$ induit un homéomorphisme $f_{\alpha'}:G(\alpha')\to G(\alpha)$. On pose alors $\mu_{\alpha'}=(f_{\alpha'}^{-1})_*\mu_\alpha$. De même, si $\alpha'$ est
 un arc 
 transverse à $\Lambda'$, alors son image est la réunion d'images de relevés d'arcs transverses à $\Lambda$, soit 
 $\alpha'([0,1])=\alpha_1'([0,1])\cup\dots\cup\alpha_n'([0,1])$ tels 
 que, pour tout $2\leqslant k\leqslant n$, l'intersection $\alpha_{k-1}'([0,1])\cap\alpha_{k}'([0,1])$ est l'image d'un relevé 
 d'un arc transverse à $\Lambda$, et pour tout $p\not\in\{k,k+1,k-1\}$, l'intersection
 $\alpha_k'([0,1])\cap\alpha_{p}'([0,1])$ est vide. 
 Alors, pour tout $k\in[1,n]\cap\NN$, on a $\mu_{\alpha_k'|G(\alpha_k')\cap G(\alpha_{k+1}')}=(f_{\alpha'_k}^{-1})_*(\mu_{\alpha_k|G(p'\circ\alpha_k')\cap
 G(p'\circ\alpha_{k+1}')})
 =(f_{\alpha'_{k+1}}^{-1})_*(\mu_{\alpha_{k+1}|G(p'\circ\alpha_k')\cap G(p'\circ\alpha_{k+1}')})=\mu_{\alpha_{k+1}'|G(\alpha_k')\cap G(\alpha_{k+1}')}$. 
 Donc il existe une unique mesure $\mu_{\alpha'}$ sur
 $G(\alpha')$ telle que $\mu_{\alpha'|G(\alpha_k')}=\mu_{\alpha_k'}$ pour tout $k\in[1,n]\cap\NN$. On note $\mu'=(\mu'_{\alpha'})_{\alpha'\in\tau'}$, où $\tau'$
 est 
 l'ensemble des arcs transverses à $\Lambda'$ la famille de mesures ainsi définie. Par naturalité, $\Lambda'$ est $\Gamma_{\Sigma'}$-invariante et la famille
 $\mu'$ est invariante pour l'action par homéomorphismes de
 $\Gamma_{\Sigma'}$ définie par ${\gamma(\mu'_{\alpha'})_{\alpha'\in\T'}}=(\gamma_*\mu'_{\gamma^{-1}\alpha'})_{\alpha'\in\T'}$, 
 pour tout $\gamma\in\Gamma_{\Sigma'}$. 
 
 \blemm\label{mu'}
 La famille $\mu'$ est l'unique mesure transverse à $\Lambda'$ telle que si $\alpha'$ est le relevé d'un arc 
 $\alpha$ transverse à $\Lambda$ alors $\mu_{\alpha'}=(f_{\alpha'}^{-1})_*\mu_\alpha$. De plus l'application de
 $\M\LLLL_p(\Sigma)$ dans $\M\LLLL_p(\Sigma')$ ainsi définie est un homéomorphisme entre $\M\LLLL_p(\Sigma)$ et l'espace des laminations plates mesurées de
 $\Sigma'$
 qui sont $\Gamma_{\Sigma'}$-invariantes. 
 \elemm

 \dem Les propriétés $(1)$, $(3)$ et $(4)$ de la définition \ref{defmesuretransverse} sont clairement satisfaites par $\mu'$, et si $h:F_1'\subset 
 G(\alpha_1')\to F_2'
 \subset G(\alpha_2')$ est une holonomie de $\Lambda'$ qui relève une holonomie $h:F_1\subset G(\alpha_1)\to F_2\subset G(\alpha_2)$ de $\Lambda$,    
  où $\alpha_1$ et $\alpha_2$, sont deux arcs disjoints et transverses à $\Lambda$, 
  alors $\mu_{\alpha_2'|F_2'}=(f_{\alpha_2'}^{-1})_*(\mu_{\alpha_2|F_2})=
 (f^{-1}_{\alpha_2'})_*h_*(\mu_{\alpha_1|F_1})=h'_*(f^{-1}_{\alpha_1'})_*(\mu_{\alpha_1|F_1})=h'_*(\mu_{\alpha'_1|F'_1})$. Sinon, on écrit $h$ comme une concaténation 
 de compositions 
 d'holonomie qui    sont des relevés d'holonomies de $\Lambda$, et on montre 
 de même que $\mu_{\alpha_2'|F_2'}=h'_*(\mu_{\alpha_1'|F_1'})$. Donc $\mu'$ est une mesure transverse à $\Lambda'$. 
 L'application de $\M\LLLL_p(\Sigma)$ dans $\M\LLLL_p(\Sigma')$ ainsi définie est injective par construction.
   Si $(\Lambda',\mu')$ est une lamination plate mesurée de $(\Sigma',[q'])$ qui est $\Gamma_{\Sigma'}$-invariante, alors l'ensemble des projetés des feuilles
  de $\Lambda'$ sur $\Sigma$ est une lamination plate $\Lambda$ de $\srfce$ et si $\alpha$ est un arc transverse à $\Lambda$ et $\alpha'$ est un relevé de $\alpha$ 
 dans $\Sigma'$, alors $(f_{\alpha'})_*\mu'_{\alpha'}$ est une mesure sur $G(\alpha)$, qui ne dépend pas du choix du relevé par $\Gamma_{\Sigma'}$-invariance. Cette
 mesure
 satisfait clairement les propriétes $(1)$, $(3)$ et $(4)$ de la définition \ref{defmesuretransverse}. De plus si $h:F_1\to F_2$ est une holonomie entre deux
 boréliens
 de $G(\alpha_1)$ et $G(\alpha_2)$, elle se relève en une holonomie entre deux boréliens de $G(\alpha_1')$ et $G(\alpha_2')$, où $\alpha_1'$ et $\alpha_2'$
 sont des 
 relevés
 de $\alpha_1$ et $\alpha_2$, et on montre que $\mu_{\alpha_2|F_2}=h_*(\mu_{\alpha_1|F_1})$. Donc $(\Lambda,\mu)$ est bien une lamination plate mesurée, et
 son image par l'application
 précédente est $(\Lambda',\mu')$.
%
%
%
%
  Donc cette application est bien une bijection
entre $\M\LLLL_p(\Sigma')$ et l'ensemble des laminations plates mesurées de $(\Sigma',[q'])$ qui sont $\Gamma_{\Sigma'}$-invariantes.
Enfin, si $(\Lambda_n,\mu_n)_{n\in\NN}$ est une suite de $\M\LLLL_p(\Sigma)$ qui converge vers $(\Lambda,\mu)$, et si $(\Lambda'_n,\mu'_n)_{n\in\NN}$ et
$(\Lambda',\mu')$ 
sont leurs images dans $\M\LLLL_p(\Sigma')$, alors si $\alpha'$ est un arc transverse à $\Lambda'$ qui est le relevé d'un arc $\alpha$ transverse à $\Lambda$,
alors
  $\alpha$ est transverse
à $\Lambda_n$ pour $n$ assez grand, et puisque $p$ est un difféomorphisme local, $\alpha'$ est transverse à $\Lambda'_n$ pour $n$ assez grand. De plus,
on a $\mu'_{n,\alpha'}=(f_{\alpha'}^{-1})_*\mu_{n,\alpha}\overset{*}{\rightharpoonup}(f_{\alpha'}^{-1})_*\mu_{\alpha}=\mu'_\alpha$. Si $\alpha'$ est un arc 
transverse à $\Lambda'$
qui n'est pas le relevé d'un arc transverse à $\Lambda$, alors en le décomposant, on montre de même que pour $n$ assez grand $\Lambda'_n$ est transverse à 
$\alpha'$
et $\mu'_{n,\alpha'}\overset{*}{\rightharpoonup}\mu'_{\alpha'}$. Donc l'application $(\Lambda,\mu)\mapsto(\Lambda',\mu')$ est continue. De
même, l'application inverse est continue, donc c'est un homéomorphisme.\cqfd

\section{Laminations plates mesurées et mesures de Radon sur l'espace des géodésiques plates du revêtement universel.}

Dans cette partie, on considère toujours une surface munie d'une structure plate $\srfce$ comme dans la partie \ref{rappelsconventions} et
un de ses revêtements universels localement isométriques, soit $p:\revet\to\srfce$.
     Les géodésiques locales de $\revet$ sont des géodésiques, et si $\Ltilde$ est une lamination plate de $\revet$ et si $\alpha$ est un arc transverse à
$\Ltilde$ tel que
chacune des géodésiques de $G(\alpha)$ n'intersecte $\alpha([0,1])$ qu'en son origine, alors l'application $g_\alpha:G(\alpha)\to[G(\alpha)]$ définie par 
$g_\alpha(g)=[g]$ est un homéomorphisme.

On note $\M_{\grperevet}([\Gqr])$ l'espace des mesures de Radon sur l'espace $[\Gqr]$  muni de la topologie des géodésiques,
    qui sont $\grperevet$ et $\iota$-invariantes (avec $\iota(\ell)=\ell^-:t\mapsto\ell(-t)$), et dont les supports sont des laminations plates 
    $\grperevet$-invariantes, muni de la topologie faible-$*$.
    Soient $\nu$ un élément de $\M_{\grperevet}([\Gqr])$
dont le support est $\Ltilde$ et $\alpha$ un arc transverse à $\Ltilde$ tel que les géodésiques de $G(\alpha)$ 
n'intersectent $\alpha([0,1])$ qu'en leur origine.  Alors $\mutilde_\alpha=(g_\alpha^{-1})_*(\nu_{|[G(\alpha)]})$
est une mesure de Radon sur $G(\alpha)$
de support $g_\alpha^{-1}(\widetilde{\Lambda}\cap[G(\alpha)])$. Si $\alpha$ est un arc transverse à $\Ltilde$, mais certaines géodésiques de $G(\alpha)$ 
intersectent
$\alpha([0,1])$ en plusieurs points, alors on définit la mesure $\mutilde_\alpha$ par recollement fini comme dans la partie \ref{relevez}, ce qui est possible
car les géodésiques de $\Gqr$ sont propres. Soit $\mutilde=\mutilde_\nu$ la famille de mesures transverses ainsi définie. 

\blemm\label{nu} La famille de mesures $\mutilde$ est une mesure transverse à $\Ltilde$, et l'application $\nu\mapsto\mutilde$ ainsi définie est un homéomorphisme entre
$\M_{\grperevet}([\Gqr])$ et
l'ensemble des laminations plates mesurées de $\revet$ qui sont $\grperevet$-invariantes (muni de la topologie induite). 
\elemm

\dem Les propriétés $(1)$, $(3)$, et $(4)$ de la définition \ref{defmesuretransverse} sont clairement satisfaites par $\mutilde$, et si 
$h:F_1\to F_2$ est une holonomie de $\Ltilde$ entre deux boréliens $F_1$ de $G(\alpha_1)$ et $F_2$ de $G(\alpha_2)$, où $\alpha_1$ et $\alpha_2$ sont des arcs transverses à 
$\Ltilde$, 
        alors, par définition des holonomies, les ensembles $[F_1]$ et $[F_2]$ sont égaux.
Donc $\mutilde_{\alpha_2|F_2}=(g_{\alpha_2}^{-1})_*\nu_{|[F_2]}=h_*(g_{\alpha_1}^{-1})_{*}\nu_{|[F_1]}=h_*(\mutilde_{\alpha_1|F_1})$ et $\mutilde$ 
est invariante 
par holonomie, donc c'est une mesure transverse à $\Ltilde$, qui est $\grperevet$-invariante par naturalité. 

Sauf mention du contraire, dans la suite de la démonstration, si $\alpha$ est un arc transverse à une lamination plate $\widetilde{\Lambda}$, on suppose
que les géodésiques de $G(\alpha)$ n'intersectent $\alpha([0,1])$ qu'en leur origine, ce à quoi on peut toujours se ramener par restriction car les géodésiques
 de $G(\alpha)$ sont propres et transverses à $\alpha$. 
 Supposons que deux mesures $\nu_1$ et $\nu_2$ de $\M_{\grperevet}([\Gqr])$ définissent la même lamination plate mesurée $(\widetilde{\Lambda},\mutilde)$ par 
cette construction. Alors si $U$ est un ouvert relativement compact de $[\Gqr]$, il existe un recouvrement ouvert fini $(U_i)_{1\leqslant i\leqslant n}$ de $U$
tel que pour chacun des 
$U_i$ il existe un arc $\alpha_i$ transverse à $\Ltilde$, tel que $U_i\subset [G(\alpha_i)]$. Mais alors
$\nu_1(U_i)=\nu_2(U_i)$ pour tout $1\leqslant i\leqslant n$, et $\nu_1(U)=\nu_2(U)$. Comme la tribu borélienne est engendrée par les ouverts relativement compacts, 
on a $\nu_1=\nu_2$. Donc l'application $\nu\mapsto\mutilde$ est injective. 

Ensuite, si $(\Ltilde,\mutilde)$ est une lamination plate mesurée, et si $U$ est un ouvert de $[\Gqr]$,
il existe une suite $(\alpha_n)_{n\in\NN}$ d'arcs transverses à $\Ltilde$ telle que $U\subset\bigcup_{n\in\NN} [G(\alpha_n)]$.
On pose $[F_1]=U\cap [G(\alpha_1)]$ et pour tout $n\geqslant2$, $[F_n]=(U-\bigcup_{k<n}[F_k]\cap[G(\alpha_n)])\cap[G(\alpha_n)]$. Pour chaque $n\in\NN$,
on note 
$F_n=\{g\in G(\alpha_n)\;:\;[g]\in[F_n]\}$. Alors les ensembles $F_n$
sont des boréliens deux-à-deux disjoints. On pose $\nu(U)=\sum_{n\in\NN}\mutilde_{\alpha_n}(F_n)$. Par invariance par holonomie, 
cette définition ne dépend pas du choix des $\alpha_n$, et on définit ainsi une 
mesure borélienne sur $[\Gqr]$. En effet, soient $(U_i)_{i\in\NN}$ une suite d'ouverts deux-à-deux disjoints et $U=\bigcup_{i\in\NN}U_i$. Il existe
une suite d'arcs 
transverses correspondant à $U$ et une décomposition 
$U=\amalg_{n\in\NN}[F_n]$ comme ci-dessus. Pour tout $n\in\NN$ on a $F_n=\amalg_{i\in\NN}F_{n,i}$, où $F_{n,i}=\{g\in F_n\;:\;[g]\in U_i\}$. 
Alors $\nu(U)=\sum_{n\in\NN}\mutilde_{\alpha_n}(F_n)=\sum_{n\in\NN}\sum_{i\in\NN}\mutilde_{\alpha_n}(F_{n,i})$, et puisque tous les termes sont positifs, on a 
$\nu(U)=\sum_{i\in\NN}\sum_{n\in\NN}\mutilde_{\alpha_n}(F_{n,i})=\sum_{i\in\NN}\nu(U_i)$. Donc la fonction ainsi définie est $\sigma$-additive et se prolonge 
de manière unique
en une mesure sur la tribu borélienne. D'après les propriétés $(1)$ et $(3)$ de la définition \ref{defmesuretransverse}, la mesure $\nu$ est de support $\Ltilde$ et 
$\iota$-invariante, elle est $\grperevet$-invariante par naturalité, et par construction, la lamination plate mesurée associée à
$\nu$ par l'application précédente est $(\Ltilde,\mutilde)$. Donc l'application $\nu\mapsto\mutilde$ est bijective.

Enfin, si $(\nu_n)_{n\in\NN}$ est une suite de $\M_{\grperevet}([\Gqr])$ qui converge vers $\nu$, et si $(\Ltilde_n,\mutilde_n)_{n\in\NN}$ et $(\Ltilde,\mutilde)$
sont leurs images dans $\M\LLLL_p(\widetilde{\Sigma})$, alors si $\alpha$ est un arc transverse à $\Ltilde$, pour $n$ suffisament grand,
$\Ltilde_n$ est transverse à $\alpha$
et $\mutilde_{n,\alpha}=(g_\alpha^{-1})_*(\nu_{n|[G(\alpha)]})\overset{*}{\rightharpoonup}(g_\alpha^{-1})_*(\nu_{|[G(\alpha)]})=\mutilde_\alpha$. Si certaines géodésiques
de $G(\alpha)$ intersectent $\alpha([0,1])$ en plus d'un point, on montre aussi que $\mutilde_{n,\alpha}\overset{*}{\rightharpoonup}\mutilde_\alpha$ 
en décomposant $\alpha([0,1])$. Donc ceci
est vrai pour tout les arcs transverses à $\Ltilde$ et la suite $(\Ltilde_n,\mutilde_n)_{n\in\NN}$ converge vers $(\Ltilde,\mutilde)$. Donc l'application est continue.
De même, si $(\Ltilde_n,\mutilde_n)_{n\in\NN}$ est une suite de laminations plates mesurées de $\widetilde{\Sigma}$ qui sont $\grperevet$-invariantes,
qui converge vers 
$(\Ltilde,\mutilde)$, et si $(\nu_n)_{n\in\NN}$ et $\nu$ sont leurs images réciproques dans $\M_{\grperevet}([\Gqr])$, alors si $f$ est une fonction continue
sur $[\Gqr]$ à valeurs dans $\RR$ et à support compact, il existe une famille finie $\{\alpha_1,\dots,\alpha_p\}$ d'arcs transverses à 
$\Ltilde$ et à $\Ltilde_n$ pour $n$ assez grand, telle que $\Supp(f)\subset\bigcup_{1\leqslant i\leqslant p}[G(\alpha_i)]$, et pour tout $1\leqslant i\leqslant p$
et $n$ assez grand, on a $\nu_n(f_{|[G(\alpha_i)]})=
\mutilde_{n,\alpha_i}(f\circ g_{\alpha_i|G(\alpha_i)})\longrightarrow\mutilde_{\alpha_i}(f\circ g_{\alpha_i|G(\alpha_i)})=\nu(f_{|[G(\alpha_i)]})$. Donc 
$\nu_n\overset{*}{\rightharpoonup}\nu$ et la réciproque est continue. Donc l'application $\nu\mapsto\mutilde$ est un homéomorphisme.\cqfd

\bcoro\label{correspondance} Si $\srfce$ est une surface munie d'une structure plate comme dans la partie \ref{rappelsconventions}, et si $p:\revet\to\srfce$ 
est un revêtement
universel localement
isométrique de groupe de revêtement $\Gamma_{\widetilde{\Sigma}}$, alors $\M\LLLL_p(\Sigma)$ est homéomorphe à $\M_{\grperevet}([\Gqr])$.\cqfd
\ecoro


\section{Liens entre laminations plates mesurées et laminations hyperboliques mesurées.}\label{liensplatshyperboliques}

 Dans cette partie \ref{liensplatshyperboliques}, on note toujours $\srfce$ une surface munie d'une structure plate comme
 dans la partie \ref{rappelsconventions}, $p:\revet\to\srfce$ un revêtement universel localement
 isométrique de groupe de revêtement $\Gamma_{\widetilde{\Sigma}}$. On suppose que $\Sigma$ est compacte et $\chi(\Sigma)<0$, et on note
 $m$ une métrique hyperbolique dans la classe conforme de $[q]$ et $\widetilde{m}$
 l'unique métrique hyperbolique
 sur $\widetilde{\Sigma}$ telle que
 $p:\revetm\to\srfcem$ soit localement isométrique. Pour toute géodésique
 $g$ de $[\Gqr]$ ou $[\Gmr]$, on note
 $E(g)\in\dddp$ le couple de ses points à l'infini, et si $F$ est un ensemble de géodésiques, $E(F)=\{E(g)\;:\;g\in F\}$. 
 L'espace des laminations hyperboliques mesurées sur $\srfcem$ (voir \cite{Bonahon97}), que l'on note $\M\LLLL_h(\Sigma)$ (muni de la topologie définie
 dans \cite[p.~19]{Bonahon97}), 
 est homéomorphe à l'espace des mesures de Radon
 $\Gamma_{\widetilde{\Sigma}}$ et $\iota$-invariantes, dont le support est une lamination hyperbolique, sur $[\Gmr]$, que l'on note $\M_{\grperevet}([\Gmr])$
 (voir
 \cite[Prop.~17~p.~154]{Bonahon88}).  Dans cette partie \ref{liensplatshyperboliques}, on utilise ce fait et l'homéomorphisme entre $\M_{\grperevet}([\Gqr])$
 et $\M\LLLL_p(\Sigma)$ défini au corollaire \ref{correspondance} pour étudier le lien entre $\M\LLLL_p(\Sigma)$ et $\M\LLLL_h(\Sigma)$.   
%
%
  On note $\varphi:[\Gqr]\to[\Gmr]$ l'application qui à une géodésique de $[\Gqr]$ associe la géodésique de $[\Gmr]$ qui lui correspond 
  (voir \cite[§4.2]{Morzy1}). 
  Alors $\varphi$ est surjective et continue, et une partie fermée $F$ de $[\Gqr]$ est une lamination plate si et seulement si $\varphi(F)$ est une lamination
  hyperbolique (voir la remarque précédant \cite[Lem.~4.3]{Morzy1}). 
  De plus $\varphi$ est 
        propre. En effet, si $K$ est un compact de $[\Gmr]$, et si $(\widetilde{\ell}_n)_{n\in\NN}$ est une suite de
  $\varphi^{-1}(K)$,
  alors par définition de $\varphi$, quitte à extraire, la suite $(E(\widetilde{\ell}_n))_{n\in\NN}$ converge dans $\dddp$. Puisque $\revet$ est 
  $\delta$-hyperbolique (car $\Sigma$ est compacte et $\chi(\Sigma)<0$, voir \cite[Rem.~2.10]{Morzy1}),
  d'après \cite[Lem.~2.6]{Morzy1}, quitte
  à extraire, la suite $(\widetilde{\ell}_n)_{n\in\NN}$ converge vers une géodésique $\widetilde{\ell}$ telle que $E(\widetilde{\ell})\in E(K)$. Donc 
  $\widetilde{\ell}$
  appartient à $\varphi^{-1}(K)$ qui est donc compact. 
  
  \blemm\label{surjectionpropre}
  L'application $\varphi$ définit une application continue, surjective et propre $\varphi_*$ entre les espaces de mesures de Radon sur $[\Gqr]$ et 
  $[\Gmr]$. De plus,
  $\varphi_*\nu$ appartient à $\M_{\grperevet}([\Gmr])$ si et seulement si $\nu$ appartient à $\M_{\grperevet}([\Gqr])$ et la restriction de $\varphi_*$ à
  $\M_{\grperevet}([\Gqr])$ est une surjection $\varphi_*:\M_{\grperevet}([\Gqr])\to\M_{\grperevet}([\Gmr])$.
  \elemm
  
 \dem L'application $\varphi$ est continue et propre donc elle définit une application continue $\varphi_*$ entre les espaces des mesures de Radon sur $[\Gqr]$ 
  et $[\Gmr]$.
%
%
       Montrons que $\varphi_*$ est surjective. L'application $s:[\Gmr]\to[\Gqr]$ qui à une géodésique hyperbolique associe la géodésique plate à laquelle elle 
  correspond
  (si elle est unique) et la géodésique "milieu" de l'ensemble des géodésiques (contenues dans une même bande plate) auxquelles elle correspond sinon,
  est une section mesurable (mais non continue) de $\varphi$.
   Puisque $\varphi$ est continue, 
  l'image réciproque d'un compact par l'application $s$ est relativement compacte, donc $s$ définit une application $s_*$ de l'ensemble des mesures de Radon sur $[\Gmr]$ dans l'ensemble 
  des 
  mesures
  de Radon sur $[\Gqr]$, et $\varphi_*\circ s_*=\Id$. Donc $\varphi_*$ est surjective.
  
  Montrons que $\varphi_*$ est propre. L'espace $[\Gmr]$ muni de la topologie des géodésiques est localement compact et $\sigma$-compact, donc il est 
  dénombrable à l'infini. Puisque $\varphi$ est continue et propre, $[\Gqr]$ est aussi dénombrable à l'infini. Donc il existe une suite $(K_n)_{n\in\NN}$ de 
  compacts telle que, pour tout $n\in\NN$, $K_n$ soit contenu dans l'intérieur de $K_{n+1}$ et $\bigcup_{n\in\NN}K_n=[\Gqr]$. Si $C$ est un compact de l'espace 
  des 
  mesures de Radon sur $[\Gmr]$
  et si $K$ est un compact de $[\Gqr]$, alors l'ensemble $\{\nu(K),\nu\in(\varphi_*)^{-1}(C)\}$ est borné par le maximum de $\{\nu(\varphi(K)),\nu\in C\}$, qui
  est fini 
  car $C$ est compact.
  Donc, pour tout $n\in\NN$, l'ensemble $\{\nu_{|K_n},\nu\in(\varphi_*)^{-1}(C)\}$ est compact. 
  Donc si $(\nu_k)_{k\in\NN}$ est 
  une
  suite de $(\varphi_*)^{-1}(C)$, par un 
  procédé d'extraction diagonale, il existe une sous-suite toujours notée $(\nu_k)_{k\in\NN}$ et une mesure de Radon $\nu$ sur $[\Gqr]$ telles que pour tout
  $n\in\NN$,
  $\nu_{k|K_n}\overset{*}{\rightharpoonup}\nu_{ |K_n}$. Or, d'après le choix de la suite $(K_n)_{n\in\NN}$, pour tout $f\in\C_c([\Gqr])$, il existe $n\in\NN$ 
  tel que 
  $\Supp(f)\subset K_n$, et alors
  $(\nu_k(f))_{k\in\NN}=(\nu_{k|K_n}(f))_{k\in\NN}$ converge vers $\nu_{|K_n}(f)=\nu(f)$. Donc on a bien $\nu_k\overset{*}{\rightharpoonup}\nu$ et $\varphi_*$  
  est propre sur les mesures de Radon.

  Enfin, par définition de $\varphi$, la mesure $\varphi_*\nu$    
           appartient à $\M_{\grperevet}([\Gmr])$ si et seulement si $\nu$ appartient à $\M_{\grperevet}([\Gqr])$. 
     De plus, l'espace $\M_{\grperevet}([\Gmr])$ est fermé (voir \cite[Prop.~3~et~17]{Bonahon88}), et puisque $\varphi_*$ est continue, 
  son image réciproque est fermée. Donc la 
  restriction de $\varphi_*$ à ces espaces définit une surjection
  continue et propre.\cqfd
  
  \medskip
  
 Les espaces $\M\LLLL_h(\Sigma)$ et $\M\LLLL_p(\Sigma)$ sont respectivement homéomorphes
  à $\M_{\grperevet}([\Gmr])$ et à $\M_{\grperevet}([\Gqr])$, donc $\varphi_*$ 
  définit une application continue, surjective et propre
  $\psi:\M\LLLL_p(\Sigma)\to\M\LLLL_h(\Sigma)$. Or, $\RR^{+*}$ agit par multiplications des mesures sur chacun des deux espaces. On note 
  $\PPPP\M\LLLL_p(\Sigma)$ et $\PPPP\M\LLLL_h(\Sigma)$
  les espaces quotients pour ces actions. Alors $\psi$ passe au quotient et définit une application continue, surjective et propre $\overline{\psi}:
  \PPPP\M\LLLL_p(\Sigma)\to
  \PPPP\M\LLLL_h(\Sigma)$. On en déduit les lemmes suivants.
  
  \blemm L'espace $\PPPP\M\LLLL_p(\Sigma)$ est compact.
  \elemm
  
  \dem L'espace $\PPPP\M\LLLL_h(\Sigma)$ est compact (voir \cite[Cor.~5~et~Prop.~17]{Bonahon88}) et $\overline{\psi}$ est propre.\cqfd
  
  \medskip
  
  Si $\Sigma$ est compacte, on appelle {\it lamination cylindrique mesurée} une lamination plate mesurée qui a une unique composante qui est cylindrique
  (voir \cite[§6]{Morzy1}).
  
  \blemm Comme $\Sigma$ est compacte, les laminations cylindriques mesurées ayant un nombre fini de feuilles sont denses dans $\M\LLLL_p(\Sigma)$. En particulier,
  $\M\LLLL_p(\Sigma)$ est séparable.
  \elemm
  
  \dem L'ensemble des géodésiques simples munies de mesures transverses, qui sont des masses de Dirac strictement positives, est dense dans $\M\LLLL_h(\Sigma)$ (voir 
  \cite[Prop.~15]{Bonahon97}), et son image réciproque par $\varphi_*$ est l'ensemble des laminations cylindriques mesurées. Puisque $\varphi_*$ est continue, 
  cet ensemble
  est dense dans $\M\LLLL_p(\Sigma)$. Si $(\Lambda,\mu)$ est une lamination cylindrique mesurée dont le support n'est pas réduit à une feuille, on note 
  $\alpha:[0,T]\to C$ ($T>0$) un arc géodésique qui relie orthogonalement les bords du cylindre plat maximal $C$ qui contient le support de $\Lambda$.
  Alors l'ensemble des géodésiques locales contenues dans $C$ et parallèles aux bords de $C$ est homéomorphe à 
  $[0,T]$, et puisque l'ensemble des mesures de Radon de support fini sur $[0,T]$ est dense dans l'espace des mesures de Radon sur $[0,T]$, il existe
  une suite $(\mu_{n,\alpha})_{n\in\NN}$ de mesures de Radon de support fini sur $G(\alpha)$ telle que pour tout $n$, chacune des feuilles du support 
  de $\mu_{n,\alpha}$ est parallèle aux bords de $C$ et $\mu_{n,\alpha}\overset{*}{\rightharpoonup}\mu_\alpha$.
  De plus, $\Lambda_n=[\Supp(\mu_{n,\alpha})]$ est une lamination plate et $\mu_{n,\alpha}$ définit une mesure transverse à $\Lambda_n$, 
  telle que la suite $(\Lambda_n,\mu_n)_{n\in\NN}$ converge vers $(\Lambda,\mu)$.\cqfd

  
\section{Nombre d'intersection.}\label{nombreintersection}
  
Dans cette partie \ref{nombreintersection}, on note toujours $\srfce$ une surface munie d'une structure plate comme dans la partie \ref{rappelsconventions} et 
on suppose que $\Sigma$ est compacte et $\chi(\Sigma)<0$. 
Un {\it courant géodésique} sur $\widetilde{\Sigma}$ est une mesure de Radon $\Gamma_{\widetilde{\Sigma}}$ et $\iota$-invariante sur $[\Gmr]$
(voir \cite[§1]{Bonahon88}). 
On note $\C(\widetilde{\Sigma})$ l'espace des courants géodésiques sur $\widetilde{\Sigma}$, muni de la topologie faible-$*$. Alors l'ensemble des classes d'homotopie
libre de courbes
fermées sur $\Sigma$ munies de masses strictement positives se plonge dans $\C(\widetilde{\Sigma})$, et le nombre d'intersection géométrique sur cet ensemble
se prolonge de manière unique en une application continue $i:\C(\widetilde{\Sigma})\times\C(\widetilde{\Sigma})\to\RR^+$ (voir \cite[Prop.~3]{Bonahon88}).
D'après le lemme \ref{surjectionpropre}, $\varphi_*$ définit une application $\varphi_*:\M_{\grperevet}([\Gqr])\to\M_{\grperevet}([\Gmr])$.   
                Puisque $\M_{\grperevet}([\Gmr])$ est homéomorphe au sous-ensemble fermé des courants géodésiques $\nu$ tels que $i(\nu,\nu)=0$ 
                (voir \cite[Prop.~17]{Bonahon88}),
 $\varphi_*$ définit une application de $\M\LLLL_p(\Sigma)$ dans l'ensemble $\{\nu\in\C(\widetilde{\Sigma})\;:\;i(\nu,\nu)=0\}$.
   Soient $\alpha$ une classe d'homotopie libre de courbes fermées, non
 triviale, $(\Lq,\mu_{[q]})$ une lamination plate mesurée et $\nu_{\mu_{[q]}}$ la mesure de $\M_{\grperevet}([\Gqr])$ 
        définie par $(\Lq,\mu_{[q]})$ (voir le corollaire \ref{correspondance}). 
 On définit le nombre d'intersection géométrique de $(\Lq,\mu_{[q]})$ et $\alpha$ par $$i_{[q]}(\mu_{[q]},\alpha)=i(\varphi_*\nu_{\mu_{[q]}},\alpha)$$
 Si $\alpha_0$ est une classe d'homotopie 
 libre
 de courbes fermées telle que 
 $\alpha=\alpha_0^k$, avec $k\in\NN$, on a $i_{[q]}(\mu_{[q]},\alpha)=k i_{[q]}(\mu_{[q]},\alpha_0)$. On suppose donc que $\alpha$ est primitive (i.e. s'il existe une classe
 d'homotopie
 libre $\alpha_0$ telle que $\alpha=\alpha_0^k$, 
 alors $k=\pm1$), on note $\alpha_{[q]}$ une géodésique locale 
 plate dans la classe de $\alpha$, $\widetilde{\alpha}_{[q]}$ un relevé de $\alpha_{[q]}$ dans $\widetilde{\Sigma}$ et $\gamma\in\grperevet-\{e\}$ 
 l'élément hyperbolique primitif de $\grperevet$ dont $\widetilde{\alpha}_{[q]}$ est un axe de translation. 
 
 \blemm\label{masseintersection} Le nombre $i_{[q]}(\mu_{[q]},\alpha)$ est égal à la moitié de la masse, pour $\nu_{\mu_{[q]}}$, de l'ensemble des feuilles de
 $\widetilde{\Lambda}$ entrelacées avec $\widetilde{\alpha}_{[q]}$, 
 qui l'intersectent
 en (au moins) un point d'un intervalle de $\widetilde{\alpha}_{[q]}(\RR)$ qui est un domaine fondamental pour l'action par translations de $\gamma^\ZZ$.
 \elemm
 
 \dem Le nombre $i(\varphi_*\nu_{\mu_{[q]}},\alpha)$ est égal à $\frac{1}{2}\varphi_*\nu_{\mu_{[q]}}(F_m)$, où $F_m$ est l'ensemble des feuilles de la 
 lamination hyperbolique mesurée
 $(\Lmr,\widetilde{\mu}_m)$ définie par $\varphi_*\nu_{\mu_{[q]}}$ qui 
 intersectent transversalement un intervalle $I=[a,\gamma a[$, avec $a\in\widetilde{\alpha}_m(\RR)$, de l'axe de translation 
 $\widetilde{\alpha}_m$ de
 $\gamma$ dans $\revetm$, qui est un domaine fondamental de $\widetilde{\alpha}_m(\RR)$ pour l'action par translations de $\gamma^\ZZ$ (voir \cite[Prop.~3]{Bonahon88}).
 Puisque $a$ est arbitraire, si $F_m$ n'est pas vide, on peut supposer que $a$ est le point d'intersection d'une feuille $\widetilde{\lambda}$ de $\Lmr$ et de
 $\widetilde{\alpha}_m$. Alors $F_m$ est l'ensemble des feuilles de $\Lmr$ qui intersectent transversalement $\widetilde{\alpha}_m$ et qui sont contenues dans 
 l'intersection des adhérences des composantes connexes, de $\widetilde{\Sigma}-\widetilde{\lambda}(\RR)$ contenant $\gamma\widetilde{\lambda}$ et de 
 $\widetilde{\Sigma}-\gamma\widetilde{\lambda}(\RR)$ contenant $\widetilde{\lambda}$, privé de $\gamma\widetilde{\lambda}$. 
 Mais $\varphi_*\nu_{\mu_{[q]}}(F_m)=\nu_{\mu_{[q]}}(F_{[q]})$, où $F_{[q]}=\varphi^{-1}(F_m)\cap\Lqr$, où $\Lqr$ est l'image réciproque de $\Lq$ dans 
 $\widetilde{\Sigma}$. Donc, par définition 
 de $\varphi$, c'est l'ensemble des feuilles de $\Lqr$ entrelacées avec $\widetilde{\alpha}_{[q]}$, qui sont contenues dans
 l'intersection des adhérences des composantes connexes de $\widetilde{\Sigma}-\widetilde{\ell}(\RR)$ contenant $\gamma\widetilde{\ell}$ et de 
 $\widetilde{\Sigma}-\gamma\widetilde{\ell}(\RR)$ contenant $\widetilde{\ell}$, privé de $\gamma\widetilde{\ell}$, 
 où $\widetilde{\ell}$ est une feuille de $\Lqr$ 
 appartenant à $\varphi^{-1}(\lambda)$ (puisque $\nu_{\mu_{[q]}}$ est $\gamma$-invariante, si plusieurs feuilles appartiennent à $\varphi^{-1}(\lambda)$, on peut choisir 
 $\widetilde{\ell}$ arbitrairement dans cet ensemble). Alors l'ensemble $F_q$ est l'ensemble des feuilles entrelacées avec $\widetilde{\alpha}_{[q]}$ qui l'intersectent 
 en au moins un point de $[b,\gamma b[$ ou $b$ est un point d'intersection de $\widetilde{\ell}$ et $\widetilde{\alpha}_{[q]}$. L'intervalle $[b,\gamma b[$
 est un domaine fondamental de $\widetilde{\alpha}_{[q]}(\RR)$ pour l'action par translations de $\gamma^\ZZ$, et puisque les choix de $\widetilde{\lambda}$ 
 et de $\widetilde{\ell}\in
 \varphi^{-1}(\lambda)$ sont arbitraires, on montre que n'importe qu'elle feuille $\widetilde{\ell}$ convient, et par $\gamma$-invariance de $\nu_{\mu_{[q]}}$,
 on montre que l'on 
 peut choisir $b$ arbitrairement dans $\widetilde{\alpha}_{[q]}(\RR)$.
 Enfin, si $F_m$ est vide, aucune feuille de $\Lqr$ n'est entrelacée avec $\widetilde{\alpha}_{[q]}$ et $i(\varphi_*\nu_{\mu_{[q]}},\alpha)=0$.\cqfd
 
 \medskip

 \rem On pourrait définir le nombre d'intersection géométrique d'une classe d'homotopie libre de courbes fermées avec $(\Lambda_{[q]},\mu_{[q]})$
 comme la borne inférieure des masses déposées par la
 lamination
 plate mesurée sur les courbes fermées, transverses à $\Lambda_{[q]}$ par morceaux, 
 dans la classe d'homotopie libre, comme pour les feuilletages 
 mesurés, 
 mais alors elle ne serait pas toujours atteinte car une géodésique locale périodique n'est généralement pas transverse par morceaux à la lamination. 
 À
 moins de définir la masse déposée par la lamination mesurée sur un segment parcouru par morceaux par des feuilles, mais cette définition est lourde et il faudrait
 distinguer beaucoup de cas.
 
 Par ailleurs, pour les surfaces munies de structures plates compactes et sans bord, à la différence des laminations hyperboliques mesurées (voir 
 \cite[Th.~2]{Ota90}),
 le nombre d'intersection avec les classes d'homotopies libres de courbes fermées de $\Sigma$ ne sépare pas les laminations plates mesurées, mais seulement leur image dans 
 $\M\LLLL_h(\Sigma)$. En particulier, la topologie définie après la définition \ref{defmesuretransverse} n'est pas
 équivalente
 à celle induite par la topologie produit sur $\RR^\HHHH$, où $\HHHH$ est l'ensemble des classes d'homotopie libre de courbes fermées, sur 
 l'image de $\M\LLLL_p(\Sigma)$  par
 l'application $(\Lambda_{[q]},\mu_{[q]})\mapsto (i(\mu_{[q]},\alpha))_{\alpha\in\HHHH}$. 
\section{Arbre dual à une lamination plate mesurée.}\label{arbredual}

Dans cette partie \ref{arbredual}, on note toujours $\srfce$ une surface munie d'une structure plate comme dans la partie \ref{rappelsconventions}, et 
l'on suppose que $\Sigma$ est compacte et $\chi(\Sigma)<0$,
$p:\revet\to\srfce$ un revêtement universel localement 
isométrique et $(\Lambda,\mu)$ une lamination plate mesurée de $\srfce$. On note $(\widetilde{\Lambda},\widetilde{\mu})$ son image réciproque dans $\revet$ et
$\nu_{\widetilde{\mu}}$
la mesure de Radon sur $\Gqr$ qui lui est 
associée (voir le lemme \ref{nu}). On suppose tout d'abord que $\nu_{\widetilde{\mu}}$ n'a pas d'atome. 

Si $\widetilde{\ell}$ est une feuille de $\widetilde{\Lambda}$, alors son complémentaire dans $\widetilde{\Sigma}$ a deux composantes connexes
et chaque feuille de $\widetilde{\Lambda}$, étant non entrelacée avec $\widetilde{\ell}$, 
est contenue dans l'une ou l'autre des adhérences des composantes connexes de $\widetilde{\Sigma}-\widetilde{\ell}(\RR)$. 
On dit qu'une feuille $\widetilde{\ell}$ {\it sépare} deux autres feuilles si celles-ci
ne sont pas contenues dans l'adhérence de la même composante connexe du complémentaire  de $\widetilde{\ell}(\RR)$ dans $\widetilde{\Sigma}$.
Soient $\widetilde{\ell}_0$ et 
$\widetilde{\ell}_1$ deux feuilles 
de $\widetilde{\Lambda}$, et $c$ un segment géodésique qui relie leur image (si les images ne sont pas disjointes $c$ peut être réduit à un point).
On note $C_i$ l'adhérence de la
composante connexe de $\widetilde{\Sigma}-\widetilde{\ell}_i(\RR)$ qui contient $\widetilde{\ell}_{i+1}(\RR)$ (avec $i\in\ZZ/2\ZZ$). On dit qu'une feuille de
$\widetilde{\Lambda}$ {\it intersecte non trivialement $c$} si elle intersecte  les deux composantes connexes du complémentaire de l'image de 
$c$ dans $C_0\cap C_1$, et on note 
$B(\widetilde{\ell}_0,\widetilde{\ell}_1)=B_{\widetilde{\Lambda}}(\widetilde{\ell}_0,\widetilde{\ell}_1)$ l'ensemble des feuilles de 
$\widetilde{\Lambda}$ qui sont contenues dans $C_0\cap C_1$ et intersectent non trivialement $c$.

\blemm L'ensemble compact $B(\widetilde{\ell}_0,\widetilde{\ell}_1)$ ne dépend pas du choix de $c$.
\elemm

\dem L'ensemble des feuilles contenues dans $C_0\cap C_1$ est compact et la condition d'intersecter non trivialement $c$ est fermée sur cet ensemble, donc 
$B(\widetilde{\ell}_0,\widetilde{\ell}_1)$ est compact.
Soit $c'$ un autre segment géodésique qui relie les images de $\widetilde{\ell}_0$ et $\widetilde{\ell}_1$. Comme $c$ et $c'$ 
séparent $C_0\cap C_1$ en deux composantes
connexes et puisque l'intersection de deux segments géodésiques de $\revet$ (ou d'une géodésique et d'un point) est connexe,
chacune des géodésiques qui intersecte non trivialement $c$ intersecte non trivialement $c'$ et réciproquement.\cqfd

\medskip

Si $\{\widetilde{\ell}_1,\widetilde{\ell}_2\}$ est une paire de feuilles de $\widetilde{\Lambda}$, on pose 
$\widetilde{d}_{\widetilde{\Lambda}}(\widetilde{\ell}_1,\widetilde{\ell}_2)=\nu_{\widetilde{\mu}}(B(\widetilde{\ell}_1,\widetilde{\ell}_2))$. Alors 
$\widetilde{d}_{\widetilde{\Lambda}}(\widetilde{\ell}_1,\widetilde{\ell}_2)\geqslant 0$ et
$\widetilde{d}_{\widetilde{\Lambda}}(\widetilde{\ell}_1,\widetilde{\ell}_2)=\widetilde{d}_{\widetilde{\Lambda}}(\widetilde{\ell}_2,\widetilde{\ell}_1)$. De plus, si 
$\widetilde{\ell}_1$, $\widetilde{\ell}_2$ et $\widetilde{\ell}_3$ sont trois feuilles de $\widetilde{\Lambda}$ et $c_1$, $c_2$ et $c_3$ sont des segments 
géodésiques
reliant respectivement les images de $\widetilde{\ell}_1$ et $\widetilde{\ell}_2$, $\widetilde{\ell}_2$ et $\widetilde{\ell}_3$, et $\widetilde{\ell}_1$ et 
$\widetilde{\ell}_3$ (noté $c_{1,2,3}$ sur la figure, cas $2$, $3$ et $4$), alors ou bien aucune des trois feuilles ne sépare les deux autres (cas $1$)
ou bien l'une sépare les deux autres
(cas $2$, $3$ et $4$). Dans le cas $1$, chacune des feuilles de $B(\widetilde{\ell}_1,\widetilde{\ell}_3)$ intersecte non trivialement $c_1$ ou $c_2$, donc 
$B(\widetilde{\ell}_1,\widetilde{\ell}_3)\subseteq B(\widetilde{\ell}_1,\widetilde{\ell}_2)\cup B(\widetilde{\ell}_2,\widetilde{\ell}_3)$, dans le cas $2$ on a 
$B(\widetilde{\ell}_1,\widetilde{\ell}_3)=B(\widetilde{\ell}_1,\widetilde{\ell}_2)\cup B(\widetilde{\ell}_2,\widetilde{\ell}_3)$, dans le cas $3$ on a 
$B(\widetilde{\ell}_1,\widetilde{\ell}_3)\subseteq B(\widetilde{\ell}_1,\widetilde{\ell}_2)$ et dans le cas $4$, on a 
$B(\widetilde{\ell}_1,\widetilde{\ell}_3)\subseteq B(\widetilde{\ell}_2,\widetilde{\ell}_3)$.
Dans tous les cas, on a $\widetilde{d}_{\widetilde{\Lambda}}(\widetilde{\ell}_1,\widetilde{\ell}_3)\leqslant\widetilde{d}_{\widetilde{\Lambda}}(\widetilde{\ell}_1,\widetilde{\ell}_2)
+\widetilde{d}_{\widetilde{\Lambda}}(\widetilde{\ell}_2,\widetilde{\ell}_3)$. Donc $\widetilde{d}_{\widetilde{\Lambda}}$ est une pseudo-distance sur $\widetilde{\Lambda}$.
\begin{center}
\begin{picture}(0,0)%
\includegraphics{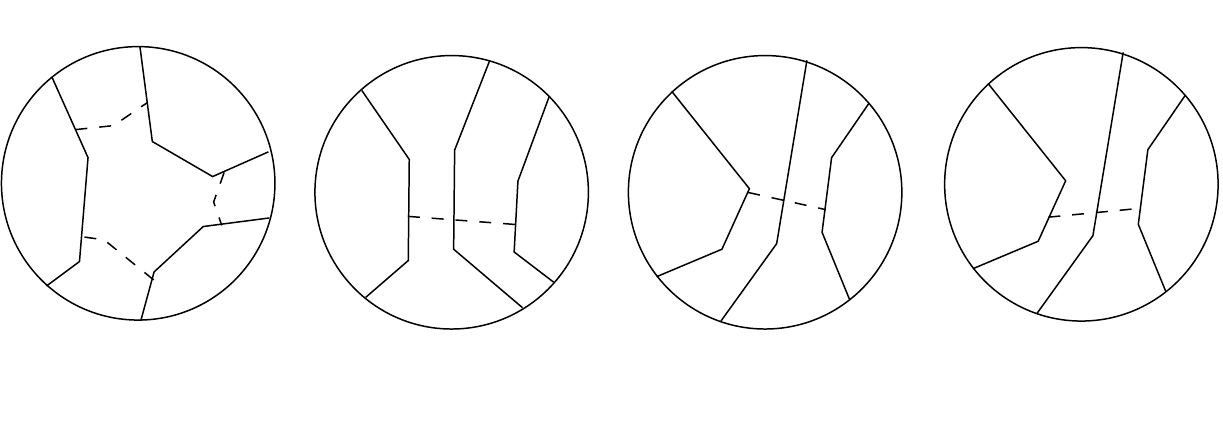}%
\end{picture}%
\setlength{\unitlength}{3771sp}%
\begingroup\makeatletter\ifx\SetFigFont\undefined%
\gdef\SetFigFont#1#2#3#4#5{%
  \reset@font\fontsize{#1}{#2pt}%
  \fontfamily{#3}\fontseries{#4}\fontshape{#5}%
  \selectfont}%
\fi\endgroup%
\begin{picture}(6129,2095)(-694,-350)
\put(951,787){\makebox(0,0)[lb]{\smash{{\SetFigFont{11}{13.2}{\familydefault}{\mddefault}{\updefault}{\color[rgb]{0,0,0}$c_{1,2,3}$}%
}}}}
\put(916,1319){\makebox(0,0)[lb]{\smash{{\SetFigFont{11}{13.2}{\familydefault}{\mddefault}{\updefault}{\color[rgb]{0,0,0}$\widetilde{\ell}_1$}%
}}}}
\put(2137,180){\makebox(0,0)[lb]{\smash{{\SetFigFont{11}{13.2}{\familydefault}{\mddefault}{\updefault}{\color[rgb]{0,0,0}$\widetilde{\ell}_3$}%
}}}}
\put(2468,1314){\makebox(0,0)[lb]{\smash{{\SetFigFont{11}{13.2}{\familydefault}{\mddefault}{\updefault}{\color[rgb]{0,0,0}$\widetilde{\ell}_1$}%
}}}}
\put(3673,1253){\makebox(0,0)[lb]{\smash{{\SetFigFont{11}{13.2}{\familydefault}{\mddefault}{\updefault}{\color[rgb]{0,0,0}$\widetilde{\ell}_2$}%
}}}}
\put(3278,1500){\makebox(0,0)[lb]{\smash{{\SetFigFont{11}{13.2}{\familydefault}{\mddefault}{\updefault}{\color[rgb]{0,0,0}$\widetilde{\ell}_3$}%
}}}}
\put(-71,-71){\makebox(0,0)[lb]{\smash{{\SetFigFont{11}{13.2}{\familydefault}{\mddefault}{\updefault}{\color[rgb]{0,0,0}$\widetilde{\ell}_2$}%
}}}}
\put(-123,560){\makebox(0,0)[lb]{\smash{{\SetFigFont{11}{13.2}{\familydefault}{\mddefault}{\updefault}{\color[rgb]{0,0,0}$c_1$}%
}}}}
\put(453,712){\makebox(0,0)[lb]{\smash{{\SetFigFont{11}{13.2}{\familydefault}{\mddefault}{\updefault}{\color[rgb]{0,0,0}$c_2$}%
}}}}
\put(-645,109){\makebox(0,0)[lb]{\smash{{\SetFigFont{11}{13.2}{\familydefault}{\mddefault}{\updefault}{\color[rgb]{0,0,0}$\widetilde{\ell}_1$}%
}}}}
\put(-210,1215){\makebox(0,0)[lb]{\smash{{\SetFigFont{11}{13.2}{\familydefault}{\mddefault}{\updefault}{\color[rgb]{0,0,0}$c_3$}%
}}}}
\put(-42,1586){\makebox(0,0)[lb]{\smash{{\SetFigFont{11}{13.2}{\familydefault}{\mddefault}{\updefault}{\color[rgb]{0,0,0}$\widetilde{\ell}_3$}%
}}}}
\put(2544,802){\makebox(0,0)[lb]{\smash{{\SetFigFont{11}{13.2}{\familydefault}{\mddefault}{\updefault}{\color[rgb]{0,0,0}$c_{1,2,3}$}%
}}}}
\put(1934,-24){\makebox(0,0)[lb]{\smash{{\SetFigFont{11}{13.2}{\familydefault}{\mddefault}{\updefault}{\color[rgb]{0,0,0}$\widetilde{\ell}_2$}%
}}}}
\put(4130,845){\makebox(0,0)[lb]{\smash{{\SetFigFont{11}{13.2}{\familydefault}{\mddefault}{\updefault}{\color[rgb]{0,0,0}$c_{1,2,3}$}%
}}}}
\put(5292,1318){\makebox(0,0)[lb]{\smash{{\SetFigFont{11}{13.2}{\familydefault}{\mddefault}{\updefault}{\color[rgb]{0,0,0}$\widetilde{\ell}_3$}%
}}}}
\put(4005,235){\makebox(0,0)[lb]{\smash{{\SetFigFont{11}{13.2}{\familydefault}{\mddefault}{\updefault}{\color[rgb]{0,0,0}$\widetilde{\ell}_2$}%
}}}}
\put(4905,1553){\makebox(0,0)[lb]{\smash{{\SetFigFont{11}{13.2}{\familydefault}{\mddefault}{\updefault}{\color[rgb]{0,0,0}$\widetilde{\ell}_1$}%
}}}}
\put(4591,-286){\makebox(0,0)[lb]{\smash{{\SetFigFont{11}{13.2}{\familydefault}{\mddefault}{\updefault}{\color[rgb]{0,0,0}cas $4$}%
}}}}
\put(-224,-286){\makebox(0,0)[lb]{\smash{{\SetFigFont{11}{13.2}{\familydefault}{\mddefault}{\updefault}{\color[rgb]{0,0,0}cas $1$}%
}}}}
\put(1351,-286){\makebox(0,0)[lb]{\smash{{\SetFigFont{11}{13.2}{\familydefault}{\mddefault}{\updefault}{\color[rgb]{0,0,0}cas $2$}%
}}}}
\put(2971,-286){\makebox(0,0)[lb]{\smash{{\SetFigFont{11}{13.2}{\familydefault}{\mddefault}{\updefault}{\color[rgb]{0,0,0}cas $3$}%
}}}}
\end{picture}%

\end{center}

On note $(T,d_T)$ l'espace métrique quotient $(\widetilde{\Lambda},\widetilde{d}_{\widetilde{\Lambda}})/\sim$ pour la relation d'équivalence
$\widetilde{\ell}\sim\widetilde{\ell}'$
si et seulement si $\widetilde{d}_{\widetilde{\Lambda}}(\widetilde{\ell},\widetilde{\ell}')=0$, et si
$F$ est un ensemble de feuilles de $\widetilde{\Lambda}$, on note $F^T$
son image par l'application de passage
au quotient. 

\brema\label{B(ell1ell2)} {\rm Soient $\widetilde{\ell}_1$ et $\widetilde{\ell}_2$ deux feuilles distinctes de $\widetilde{\Lambda}$. 
Comme pour les images réciproques des laminations hyperboliques mesurées sur une surface compacte, si $\widetilde{\ell}_1$ et
$\widetilde{\ell}_2$ ne bordent pas une même bande plate et sont telles que $\widetilde{\ell}_1(+\infty)=\widetilde{\ell}_2(+\infty)$, alors il n'existe pas d'autre 
feuille $\widetilde{\ell}$ de $\widetilde{\Lambda}$ telle que $\widetilde{\ell}(+\infty)=\widetilde{\ell}_1(+\infty)$. On a alors 
$B(\widetilde{\ell}_1,\widetilde{\ell}_2)=\{\widetilde{\ell}_1,\widetilde{\ell}_2\}$.     
 
 De plus, comme $\nu_{\widetilde{\mu}}$ est de support $\widetilde{\Lambda}$ et sans atome, si $\widetilde{\ell}_1$ et $\widetilde{\ell}_2$ bordent une même bande plate
 d'intérieur non vide ou ont
 leurs points à l'infini qui sont deux-à-deux distincts, alors 
 $B(\widetilde{\ell}_1,\widetilde{\ell}_2)$ est ou bien réduit à 
 $\{\widetilde{\ell}_1,\widetilde{\ell}_2\}$, ou bien d'intérieur non vide. En effet, si $\widetilde{\ell}$ est une feuille de 
 $B(\widetilde{\ell}_1,\widetilde{\ell}_2)-\{\widetilde{\ell}_1,\widetilde{\ell}_2\}$, alors il existe un voisinage ouvert $U$ de $\widetilde{\ell}$ dans 
 $\widetilde{\Lambda}$, non réduit à $\widetilde{\ell}$, et par hypothèse, quitte à restreindre, $U$ est contenu dans l'intérieur de 
 ${B}(\widetilde{\ell}_1,\widetilde{\ell}_2)$.
 
 Donc, dans tous les cas, ou bien $B(\widetilde{\ell}_1,\widetilde{\ell}_2)=\{\widetilde{\ell}_1,\widetilde{\ell}_2\}$,
 ou bien $B(\widetilde{\ell}_1,\widetilde{\ell}_2)$ est d'intérieur
 non vide (dans $\widetilde{\Lambda}$).
 Puisque $\nu_{\mutilde}$ n'a pas d'atome et son support est $\widetilde{\Lambda}$, on en conclut que $\nu_{\mutilde}(B(\widetilde{\ell}_1,\widetilde{\ell}_2))=0$ si
 et seulement si $B(\widetilde{\ell}_1,\widetilde{\ell}_2)=\{\widetilde{\ell}_1,\widetilde{\ell}_2\}$, et que la topologie définie par la distance $d_T$ est équivalente à la
topologie quotient de la topologie des géodésiques sur $\widetilde{\Lambda}$, pour la relation $\widetilde{\ell}_1\;R\;\widetilde{\ell}_2$ si et seulement si 
$B(\widetilde{\ell}_1,\widetilde{\ell}_2)=\{\widetilde{\ell}_1,\widetilde{\ell}_2\}$.

Enfin, puisque le complémentaire de chacune des feuilles a deux
composantes
connexes et chacune des autres feuilles est contenue dans l'adhérence de l'une des deux, la relation définie sur $B(\widetilde{\ell}_1,\widetilde{\ell}_2)$ par 
$\widetilde{\ell}\preceq\widetilde{\ell}'$ si et seulement si $\widetilde{\ell}$ appartient à $B(\widetilde{\ell}_1,\widetilde{\ell}')$ est une relation 
d'ordre total qui passe au quotient et définit une relation d'ordre total sur $B(\widetilde{\ell}_1,\widetilde{\ell}_2)^T$.}

\erema

\blemm\label{arbredualaunelamination} L'espace métrique $(T,d_T)$ est un arbre réel.
\elemm

\dem Soient $\widetilde{\ell}_1$ et $\widetilde{\ell}_2$ deux feuilles de $\widetilde{\Lambda}$.
       L'application
$f:B(\widetilde{\ell}_1,\widetilde{\ell}_2)\to\RR^+$ 
définie par $f(\widetilde{\ell})=\widetilde{d}_{\widetilde{\Lambda}}(\widetilde{\ell}_1,\widetilde{\ell})$ est croissante (pour $\preceq$) et continue car
$\nu_{\widetilde{\mu}}$ n'a pas d'atome. 
De plus elle passe au quotient et définit une application continue et strictement croissante $\overline{f}:B(\widetilde{\ell}_1,\widetilde{\ell}_2)^T\to\RR^+$.
Comme $B(\widetilde{\ell}_1,\widetilde{\ell}_2)^T$ est compact, c'est un homéomorphisme sur son image. Supposons que son image ne soit pas un intervalle. Puisque
c'est un sous-ensemble compact de $\RR$, si $U$ est une composante connexe bornée du complémentaire de ${f}(B(\widetilde{\ell}_1,\widetilde{\ell}_2))$
dans $\RR^+$, alors son adhérence est un intervalle $[a,b]$ avec $a<b$. Soient $\widetilde{\ell}_a$ et $\widetilde{\ell}_b$ des feuilles de
${f}^{-1}(a)$ et ${f}^{-1}(b)$. Si $B(\widetilde{\ell}_a,\widetilde{\ell}_b)-\{\widetilde{\ell}_a,\widetilde{\ell}_b\}$ n'est pas vide, alors on a vu que 
$B(\widetilde{\ell}_a,\widetilde{\ell}_b)$ est d'intérieur non vide. Mais si $\widetilde{\ell}$ appartient à l'intérieur de 
$B(\widetilde{\ell}_a,\widetilde{\ell}_b)$, les ensembles $B(\widetilde{\ell}_1,\widetilde{\ell})$ et $B(\widetilde{\ell},\widetilde{\ell}_2)$ sont d'intérieur non vide,
et $a<\overline{f}(\widetilde{\ell}^T)<b$. Donc  
$B(\widetilde{\ell}_a,\widetilde{\ell}_b)=\{\widetilde{\ell}_a,\widetilde{\ell}_b\}$, ce qui est impossible
car alors $\widetilde{d}(\widetilde{\ell}_{a},\widetilde{\ell}_{b})=0$ et on aurait $a=b$. Donc l'image de $\overline{f}$ est l'intervalle 
$[0,d_T(\widetilde{\ell}^T$,$\widetilde{\ell}'^T)]$,
et $\overline{f}^{-1}:[0,d_T(\widetilde{\ell}^T$, $\widetilde{\ell}'^T)]\to T$ est un segment entre $\widetilde{\ell}^T$ et 
$\widetilde{\ell}'^T$. 
Montrons que c'est le seul arc entre $\widetilde{\ell}^T$ et $\widetilde{\ell}'^T$, à reparamétrage près. Soit $g:[0,1]\to T$ un autre arc 
joignant $\widetilde{\ell}^T$ et $\widetilde{\ell}'^T$. 
Si une feuille $\widetilde{\ell}_0$ appartient à $B(\widetilde{\ell},\widetilde{\ell}')$ alors elle sépare $\widetilde{\ell}$ et $\widetilde{\ell}'$ 
(au sens de la définition ci-dessus) et puisque
$g$ est continue pour la topologie quotient de la topologie des géodésiques par la relation $R$, $\widetilde{\ell}_0^T$
appartient à l'image de $g$. Donc $B(\widetilde{\ell},\widetilde{\ell}')^T$ est contenu dans l'image de $g$. Supposons qu'il existe un élément $x$ 
de l'image de $g$ 
qui n'appartienne pas à
$B(\widetilde{\ell},\widetilde{\ell}')^T$ et soit $\widetilde{\ell}_x$ une feuille représentant $x$. Si $\widetilde{\ell}_x$ était contenue dans 
l'intersection des
adhérences 
des composantes connexes des complémentaires des images de $\widetilde{\ell}$ et $\widetilde{\ell}'$ qui contiennent respectivement $\widetilde{\ell}'$ et 
$\widetilde{\ell}$ (cas $1$ ci-dessous), alors puisque $g$ est continue, il existerait un élément de l'image de $g$ dont un représentant serait entrelacé avec 
$\widetilde{\ell}$ ou $\widetilde{\ell}'$, ce qui est impossible car les feuilles de $\widetilde{\Lambda}$ sont deux-à-deux non entrelacées. 
Donc $\widetilde{\ell}_x$ est contenue dans l'adhérence de la 
composante connexe du complémentaire de $\widetilde{\ell}$ ne contenant pas $\widetilde{\ell}'$ ou l'inverse (cas $2$ et $3$ ci-dessous). 
Mais puisque $g$ est continue et comme $\widetilde{\ell}$ sépare $\widetilde{\ell}_x$ de $\widetilde{\ell}'$ (ou $\widetilde{\ell}'$ de $\widetilde{\ell}$),
il existerait alors $t\in\mathopen{]}0,1[$ tel que $\widetilde{\ell}$ (ou $\widetilde{\ell}'$) représente $g(t)$, et $g$ ne serait pas
injectif.
\begin{center}
 \begin{picture}(0,0)%
\includegraphics{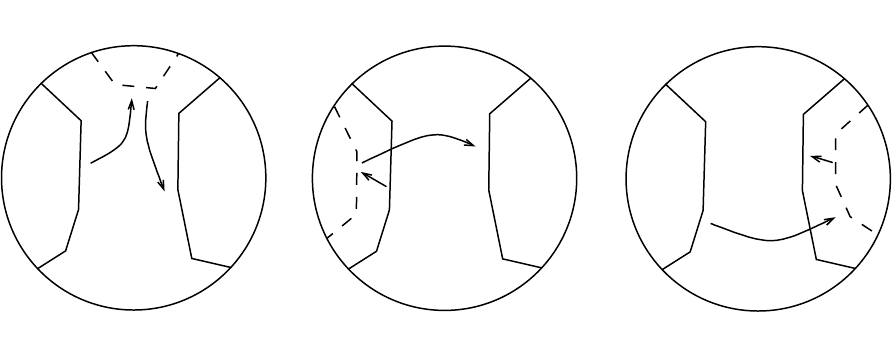}%
\end{picture}%
\setlength{\unitlength}{3771sp}%
\begingroup\makeatletter\ifx\SetFigFont\undefined%
\gdef\SetFigFont#1#2#3#4#5{%
  \reset@font\fontsize{#1}{#2pt}%
  \fontfamily{#3}\fontseries{#4}\fontshape{#5}%
  \selectfont}%
\fi\endgroup%
\begin{picture}(4482,1799)(-671,1719)
\put(991,1964){\makebox(0,0)[lb]{\smash{{\SetFigFont{11}{13.2}{\familydefault}{\mddefault}{\updefault}{\color[rgb]{0,0,0}$\widetilde{\ell}$}%
}}}}
\put(-314,3359){\makebox(0,0)[lb]{\smash{{\SetFigFont{11}{13.2}{\familydefault}{\mddefault}{\updefault}{\color[rgb]{0,0,0}$\widetilde{\ell}_x$}%
}}}}
\put(451,3179){\makebox(0,0)[lb]{\smash{{\SetFigFont{11}{13.2}{\familydefault}{\mddefault}{\updefault}{\color[rgb]{0,0,0}$\widetilde{\ell}'$}%
}}}}
\put(721,2279){\makebox(0,0)[lb]{\smash{{\SetFigFont{11}{13.2}{\familydefault}{\mddefault}{\updefault}{\color[rgb]{0,0,0}$\widetilde{\ell}_x$}%
}}}}
\put(2071,2009){\makebox(0,0)[lb]{\smash{{\SetFigFont{11}{13.2}{\familydefault}{\mddefault}{\updefault}{\color[rgb]{0,0,0}$\widetilde{\ell}'$}%
}}}}
\put(3601,2594){\makebox(0,0)[lb]{\smash{{\SetFigFont{11}{13.2}{\familydefault}{\mddefault}{\updefault}{\color[rgb]{0,0,0}$\widetilde{\ell}_x$}%
}}}}
\put(-629,1964){\makebox(0,0)[lb]{\smash{{\SetFigFont{11}{13.2}{\familydefault}{\mddefault}{\updefault}{\color[rgb]{0,0,0}$\widetilde{\ell}$}%
}}}}
\put(1306,1784){\makebox(0,0)[lb]{\smash{{\SetFigFont{11}{13.2}{\familydefault}{\mddefault}{\updefault}{\color[rgb]{0,0,0}cas $2$}%
}}}}
\put(2926,1784){\makebox(0,0)[lb]{\smash{{\SetFigFont{11}{13.2}{\familydefault}{\mddefault}{\updefault}{\color[rgb]{0,0,0}cas $3$}%
}}}}
\put(-179,1784){\makebox(0,0)[lb]{\smash{{\SetFigFont{11}{13.2}{\familydefault}{\mddefault}{\updefault}{\color[rgb]{0,0,0}cas $1$}%
}}}}
\put(2566,3134){\makebox(0,0)[lb]{\smash{{\SetFigFont{11}{13.2}{\familydefault}{\mddefault}{\updefault}{\color[rgb]{0,0,0}$\widetilde{\ell}$}%
}}}}
\put(3646,3134){\makebox(0,0)[lb]{\smash{{\SetFigFont{11}{13.2}{\familydefault}{\mddefault}{\updefault}{\color[rgb]{0,0,0}$\widetilde{\ell}'$}%
}}}}
\end{picture}%

\end{center}

Donc si $g$ est un arc entre $\widetilde{\ell}$ et $\widetilde{\ell}'$, il a la même image que $\overline{f}^{-1}$.
L'unique arc (à reparamétrage près) entre $\widetilde{\ell}^T$ et $\widetilde{\ell}'^T$ est donc $\overline{f}^{-1}$ qui est isométrique à 
$[0,d_T(\widetilde{\ell}^T$,$\widetilde{\ell}'^T)]$ par construction. Ceci étant vrai pour toutes les paires de feuilles, 
$(T,d_T)$ est un arbre réel.\cqfd

\medskip

Supposons que la mesure $\nu_{\widetilde{\mu}}$ associée à $(\widetilde{\Lambda},\widetilde{\mu})$ ait un atome $\widetilde{\ell}$.
       On remplace alors $\widetilde{\ell}(\RR)$ par une bande plate de largeur $\nu_{\widetilde{\mu}}(\widetilde{\ell})$ en recollant isométriquement chacune des 
adhérences des 
composantes
connexes du complémentaire de $\widetilde{\ell}(\RR)$ dans $\widetilde{\Sigma}$ sur 
les bords d'une copie, notée $BP(\widetilde{\ell})$, de $\RR\times[0,\nu_{\widetilde{\mu}}(\widetilde{\ell})]$ munie de 
la distance 
euclidienne. En procédant ainsi pour chacun des atomes de $\nu_{\widetilde{\mu}}$, on obtient une surface munie d'une structure plate $(\widetilde{\Sigma}',
[\widetilde{q}]')$ et l'action par 
isométries du groupe de revêtement $\Gamma_{\widetilde{\Sigma}}$ sur $\widetilde{\Sigma}$ privé des images des atomes de $\nu_{\widetilde{\mu}}$ se prolonge
de manière unique en une 
action par 
isométries sur $(\widetilde{\Sigma}',[\widetilde{q}]')$. 

Soient $\widetilde{\ell}$ un atome de $\nu_{\widetilde{\mu}}$, $F_{\widetilde{\ell}}$ l'ensemble maximal des géodésiques de la bande plate 
$\BP(\widetilde{\ell})$ de $(\widetilde{\Sigma}',[\widetilde{q}]')$ correspondante qui sont parallèles à ses bords, et  $\alpha$ un segment géodésique
de $\BP(\widetilde{\ell})$ qui relie orthogonalement les bords de $\BP(\widetilde{\ell})$. Alors
l'application $r:F_{\widetilde{\ell}}\to\Image(\alpha)$ définie par
$r(g)=g(\RR)\cap\Image(\alpha)$ est un homéomorphisme. Donc on peut munir $F_{\widetilde{\ell}}$ de la mesure $\nu_{\widetilde{\ell}}=
(r^{-1})_*dx_ \alpha$, où 
$dx_\alpha$ est la mesure proportionnelle à la mesure de Lebesgue sur $\Image(\alpha)$, de masse $\nu_{\widetilde{\mu}}(\widetilde{\ell})$. 
On définit alors la lamination plate mesurée $\grperevet$-invariante $(\widetilde{\Lambda}',\widetilde{\mu}')$ de $(\widetilde{\Sigma}',[\widetilde{q}]')$
en remplaçant les atomes 
$\widetilde{\ell}$ de $\nu_{\widetilde{\mu}}$ 
par les
ensembles $F_{\widetilde{\ell}}$ munis des mesures boréliennes $\nu_{\widetilde{\ell}}$. L'application de $\widetilde{\Lambda}'$ dans $\widetilde{\Lambda}$ induite par 
le plongement isométrique canonique de chacune des composantes connexes du complémentaire, dans $\widetilde{\Sigma}$,
de $\bigcup_{\nu_{\mutilde}(\widetilde{\ell})>0}\BP(\widetilde{\ell})$, dans $\widetilde{\Sigma}'$, et   
telle 
que l'image de toutes les feuilles d'une bande plate $\BP(\widetilde{\ell})$ soit égale à $\widetilde{\ell}$ (si $\widetilde{\ell}$ est un atome de
$\nu_{\mutilde}$), est continue, surjective et $\grperevet$-équivariante, et la mesure image de $\nu_{\mutilde'}$ par cette application est $\nu_{\mutilde}$. 
Alors $(\widetilde{\Lambda}',\widetilde{\mu}')$ n'a pas d'atome et on définit l'arbre dual à $(\widetilde{\Lambda},\widetilde{\mu})$ comme étant l'arbre dual à
$(\widetilde{\Lambda}',\widetilde{\mu}')$.

\section{Action du groupe de revêtement sur l'arbre dual à une lamination plate mesurée.}\label{actiondegroupe}

Dans cette partie, on reprend les définitions et notations de la partie \ref{arbredual} et on définit l'action canonique du groupe de revêtement 
$\Gamma_{\widetilde{\Sigma}}$ sur l'arbre dual $(T,d_T)$ à l'image réciproque $(\widetilde{\Lambda},\widetilde{\mu})$ de la lamination plate mesurée 
$(\Lambda,\mu)$ dans $\widetilde{\Sigma}$. On peut toujours supposer que $\nu_{\widetilde{\mu}}$ n'a pas d'atome, quitte à procéder comme au dernier
paragraphe de la partie \ref{arbredual}. Puisque $\widetilde{\Lambda}$ est fixée, 
on notera 
$B(\widetilde{\ell},\widetilde{\ell}')=B_{\widetilde{\Lambda}}(\widetilde{\ell},\widetilde{\ell}')$ pour toutes les feuilles $\widetilde{\ell}$, 
$\widetilde{\ell}'$ de 
$\widetilde{\Lambda}$. 

\medskip

Le groupe de revêtement $\Gamma_{\widetilde{\Sigma}}$ agit sur $\widetilde{\Sigma}$ par isométries, donc définit une action sur l'ensemble 
$[\G_{[\widetilde{q}]}]$
des géodésiques 
de $\revet$ définies à changements d'origines près. Par $\grperevet$-invariance, cette action définit une action sur $\widetilde{\Lambda}$, 
et comme pour tout $\gamma\in\grperevet$ on a $\gamma_*\nu_{\widetilde{\mu}}=\nu_{\widetilde{\mu}}$ et 
$\gamma B(\widetilde{\ell},\widetilde{\ell}')=B(\gamma\widetilde{\ell},\gamma\widetilde{\ell}')$, pour toutes  
les feuilles $\widetilde{\ell}$ et $\widetilde{\ell}'$ de $\widetilde{\Lambda}$, elle passe au quotient
       et définit une action par isométries de $\grperevet$ sur l'arbre dual $(T,d_T)$ à $(\widetilde{\Lambda},\widetilde{\mu})$ défini au lemme 
\ref{arbredualaunelamination}.

\blemm Pour tout $\gamma\in\grperevet-\{e\}$, si $\widetilde{\alpha}_{\gamma}$ est un axe de translation de $\gamma$ dans $\revet$ et $\alpha_\gamma$
est la projection de $[x,\gamma x]$, pour $x\in\widetilde{\alpha}_\gamma(\RR)$,  
dans $\Sigma$, alors la distance de translation $\ell_T(\gamma)$ de $\gamma$ dans $(T,d_T)$ est égale à $i_{[q]}(\mu,\alpha_\gamma)$. De plus, si 
$\ell_T(\gamma)>0$, l'axe de translation 
de $\gamma$ est l'image dans $T$ de l'ensemble des feuilles de $\widetilde{\Lambda}$ qui sont entrelacées avec $\widetilde{\alpha}_\gamma$.  
\elemm

\dem Il suffit de démontrer le lemme dans le cas où $\gamma$ est primitif. 

{\bf Cas $(1)$.~} Supposons que $i_{[q]}(\mu,\alpha_\gamma)>0$. Alors $\alpha_\gamma$ est entrelacée
avec au moins une feuille de 
$\Lambda$. Soit $I=[a,\gamma a[$ un domaine fondamental de $\widetilde{\alpha}_\gamma(\RR)$ pour l'action
par translations de $\gamma^\ZZ$
tel que $a$ appartienne à une feuille $\widetilde{\ell}$ de $\widetilde{\Lambda}$ entrelacée avec $\widetilde{\alpha}_\gamma$. 
Alors l'ensemble $F$ des feuilles de $\widetilde{\Lambda}$ contenues dans l'intersection des adhérences de la composante connexe de 
$\widetilde{\Sigma}-\widetilde{\ell}(\RR)$ contenant $\gamma\widetilde{\ell}$ et de la composante connexe de 
$\widetilde{\Sigma}-\gamma\widetilde{\ell}(\RR)$ contenant $\widetilde{\ell}$, privé de $\gamma\widetilde{\ell}$, est un domaine fondamental de
$\widetilde{\Lambda}$ pour l'action de $\gamma^\ZZ$. De plus, le sous-ensemble des feuilles de $F$ qui sont entrelacées avec $\widetilde{\alpha}_\gamma$ 
(et donc intersectent non trivialement $[a,\gamma a]$)
est égal à 
$B(\widetilde{\ell},\gamma\widetilde{\ell})-\gamma\widetilde{\ell}$, donc d'après le lemme \ref{masseintersection} et 
comme $\nu_{\widetilde{\mu}}(\{\gamma\widetilde{\ell}\})=0$, on a $i_{[q]}(\alpha_\gamma,\mu)=\nu_{\mutilde}(B(\widetilde{\ell},\gamma\widetilde{\ell}))
=\widetilde{d}(\widetilde{\ell},\gamma\widetilde{\ell})$.

\noindent
 \begin{minipage}{10 cm}
Comme $F$ est un ensemble fondamental de $\widetilde{\Lambda}$ pour l'action de $\gamma^\ZZ$, si $\widetilde{\ell}'$ est une feuille de 
$\widetilde{\Lambda}$, il existe un unique $n\in\ZZ$ tel que $\gamma^n\widetilde{\ell}'$ appartienne à $F$. 
On note $\widetilde{\ell}'_0=\gamma^n\widetilde{\ell}'$, $\widetilde{\ell}_1'=\gamma\widetilde{\ell}'_0$ et 
$\widetilde{\ell}'_{-1}=\gamma^{-1}\widetilde{\ell}'_0$. Alors $\widetilde{\ell}_{-1}'$ et $\widetilde{\ell}_1'$ n'appartiennent pas à $F$. Comme 
les feuilles de $\widetilde{\Lambda}$ sont deux-à-deux non entrelacées, les feuilles de $B(\widetilde{\ell},\gamma\widetilde{\ell})$
appartiennent à $B(\widetilde{\ell}'_{-1},\widetilde{\ell}'_{0})\cup B(\widetilde{\ell}'_{0},\widetilde{\ell}'_{1})$, et
$\gamma (B(\widetilde{\ell}'_{-1},\widetilde{\ell}'_0)\cap B(\widetilde{\ell},\gamma\widetilde{\ell})-\widetilde{\ell})\subseteq B(\widetilde{\ell}'_{0},
\widetilde{\ell}'_{1})-B(\widetilde{\ell}'_{0},\widetilde{\ell}'_{1})\cap B(\widetilde{\ell},\gamma\widetilde{\ell})$, 
   car $F$ est un domaine fondamental de $\widetilde{\Lambda}$ pour l'action de $\gamma^\ZZ$. Par $\grperevet$-invariance de $\nu_{\mutilde}$,
on a
\end{minipage}
\begin{minipage}{4.9 cm}
\begin{picture}(0,0)%
\includegraphics{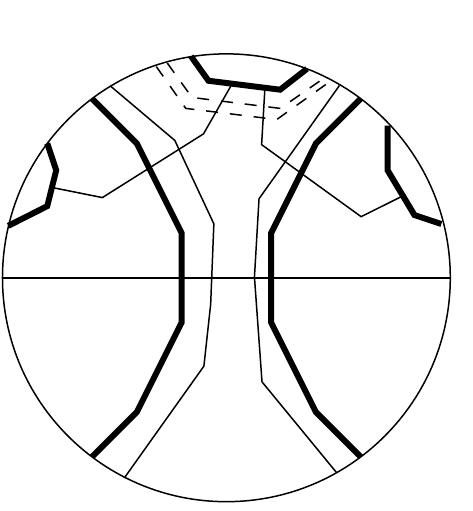}%
\end{picture}%
\setlength{\unitlength}{3771sp}%
\begingroup\makeatletter\ifx\SetFigFont\undefined%
\gdef\SetFigFont#1#2#3#4#5{%
  \reset@font\fontsize{#1}{#2pt}%
  \fontfamily{#3}\fontseries{#4}\fontshape{#5}%
  \selectfont}%
\fi\endgroup%
\begin{picture}(2362,2549)(-1136,-329)
\put( 16,2061){\makebox(0,0)[lb]{\smash{{\SetFigFont{11}{13.2}{\familydefault}{\mddefault}{\updefault}{\color[rgb]{0,0,0}$\widetilde{\ell}_{0}'$}%
}}}}
\put(-39,1457){\makebox(0,0)[lb]{\smash{{\SetFigFont{11}{13.2}{\familydefault}{\mddefault}{\updefault}{\color[rgb]{0,0,0}$A$}%
}}}}
\put(1040,1408){\makebox(0,0)[lb]{\smash{{\SetFigFont{11}{13.2}{\familydefault}{\mddefault}{\updefault}{\color[rgb]{0,0,0}$\widetilde{\ell}_{1}'$}%
}}}}
\put(675,-252){\makebox(0,0)[lb]{\smash{{\SetFigFont{11}{13.2}{\familydefault}{\mddefault}{\updefault}{\color[rgb]{0,0,0}$\gamma\cdot\widetilde{\ell}$}%
}}}}
\put(-855,-260){\makebox(0,0)[lb]{\smash{{\SetFigFont{11}{13.2}{\familydefault}{\mddefault}{\updefault}{\color[rgb]{0,0,0}$\widetilde{\ell}$}%
}}}}
\put(-404,659){\makebox(0,0)[lb]{\smash{{\SetFigFont{11}{13.2}{\familydefault}{\mddefault}{\updefault}{\color[rgb]{0,0,0}$a$}%
}}}}
\put(1211,795){\makebox(0,0)[lb]{\smash{{\SetFigFont{11}{13.2}{\familydefault}{\mddefault}{\updefault}{\color[rgb]{0,0,0}$\widetilde{\alpha}_\gamma$}%
}}}}
\put(316,659){\makebox(0,0)[lb]{\smash{{\SetFigFont{11}{13.2}{\familydefault}{\mddefault}{\updefault}{\color[rgb]{0,0,0}$\gamma a$}%
}}}}
\put(-899,929){\makebox(0,0)[lb]{\smash{{\SetFigFont{11}{13.2}{\familydefault}{\mddefault}{\updefault}{\color[rgb]{0,0,0}$\widetilde{\ell}_{-1}'$}%
}}}}
\end{picture}%

\end{minipage}
 
\begin{align*}
\nu_{\mutilde}(B(\widetilde{\ell}',\gamma\widetilde{\ell}'))&=\nu_{\mutilde}(B(\widetilde{\ell}'_{0},\widetilde{\ell}'_1))\\
&\geqslant \nu_{\mutilde}((B(\widetilde{\ell},\gamma\widetilde{\ell})\cap B(\widetilde{\ell}'_0,\widetilde{\ell}'_1))+
\nu_{\mutilde}(\gamma(B(\widetilde{\ell},\gamma\widetilde{\ell})\cap B(\widetilde{\ell}'_{-1},\widetilde{\ell}'_0)))\\
&=\nu_{\mutilde}((B(\widetilde{\ell},\gamma\widetilde{\ell})\cap B(\widetilde{\ell}'_0,\widetilde{\ell}'_1))+
\nu_{\mutilde}(B(\widetilde{\ell},\gamma\widetilde{\ell})\cap B(\widetilde{\ell}'_{-1},\widetilde{\ell}'_0))\\
&\geqslant\nu_{\mutilde}(B(\widetilde{\ell},\gamma\widetilde{\ell}))\mbox{ car } B(\widetilde{\ell},\gamma\widetilde{\ell})\subseteq 
B(\widetilde{\ell}'_{-1},\widetilde{\ell}'_0)\cup B(\widetilde{\ell}'_0,\widetilde{\ell}'_1)\\
&=\widetilde{d}(\widetilde{\ell},\gamma\widetilde{\ell})
\end{align*}

Donc $\ell_T(\gamma)=\widetilde{d}(\widetilde{\ell},\gamma\widetilde{\ell})=i_{[q]}(\mu,\alpha_\gamma)$. En particulier, on a $\ell_T(\gamma)>0$ et
l'isométrie $\gamma$ de $(T,d_T)$ est hyperbolique, et donc admet un axe de translation.   
De plus, si $\widetilde{\ell'}$ n'est pas entrelacée avec 
$\widetilde{\alpha}_\gamma$, alors les feuilles $\widetilde{\ell}'_{-1}$, $\widetilde{\ell}'_0$ et $\widetilde{\ell}_1$ non plus. Or,
on a vu que $B(\widetilde{\ell},\gamma\widetilde{\ell})$ est contenu dans $B(\widetilde{\ell}'_{-1},\widetilde{\ell}'_0)\cup B(\widetilde{\ell}'_0,\widetilde{\ell}'_1)$ et
que $ B(\widetilde{\ell}'_0,\widetilde{\ell}'_1)$ contient la réunion des ensembles $B(\widetilde{\ell},\gamma\widetilde{\ell})\cap 
B(\widetilde{\ell}'_0,\widetilde{\ell}'_1)$ et $\gamma(B(\widetilde{\ell},\gamma\widetilde{\ell})\cap 
B(\widetilde{\ell}'_{-1},\widetilde{\ell}'_0))$ dont l'intersection est  
$\{\gamma\widetilde{\ell}\}$.     
Supposons que $\widetilde{d}(\widetilde{\ell}'_0,\widetilde{\ell}'_{1})$ soit égale à $\widetilde{d}(\widetilde{\ell},\gamma\widetilde{\ell})$, i.e.
$\nu_{\mutilde}(B(\widetilde{\ell}'_0,\widetilde{\ell}'_1))=\nu_{\mutilde}(B(\widetilde{\ell},\gamma\widetilde{\ell})\cap 
B(\widetilde{\ell}'_0,\widetilde{\ell}'_1))+\nu_{\mutilde}(\gamma(B(\widetilde{\ell},\gamma\widetilde{\ell})\cap 
B(\widetilde{\ell}'_{-1},\widetilde{\ell}'_0)))$. Puisque $\nu_{\mutilde}$ est de support $\widetilde{\Lambda}$, 
cela signifie que l'intersection 
$(F\cup\gamma\widetilde{\ell})\cap(B(\widetilde{\ell}'_{0},\widetilde{\ell}_1)
\cup \gamma^{-1}B(\widetilde{\ell}'_0,\widetilde{\ell}'_1))
=(F\cup\gamma\widetilde{\ell})\cap(B(\widetilde{\ell}'_{-1},\widetilde{\ell}_0)\cup B(\widetilde{\ell}'_0,\widetilde{\ell}'_1))$ est égale
à $B(\widetilde{\ell},\gamma\widetilde{\ell})$ et donc qu'il n'existe pas de  
feuille de $F$ contenue dans l'adhérence $A$ du complémentaire de $\widetilde{\ell}'_0(\RR)$ qui contient 
$\widetilde{\alpha}_\gamma$, qui ne soit pas entrelacée avec $\widetilde{\alpha}_\gamma$ (il n'y a pas de feuille telle qu'en pointillés sur la figure).
Mais alors ou bien l'ensemble $B(\widetilde{\ell}'_0,\widetilde{\ell})$ est réduit à $\{\widetilde{\ell}'_0,\widetilde{\ell}\}$ et $\widetilde{d}(\widetilde{\ell},\widetilde{\ell}'_0)=0$, 
ou bien il a un minimum (pour la 
relation d'ordre total $\preceq$ sur $B(\widetilde{\ell}'_0,\widetilde{\ell})$ définie à la remarque \ref{B(ell1ell2)}) 
distinct de $\widetilde{\ell}'_0$, soit $m$.
Alors $\widetilde{d}(m,\widetilde{\ell}'_0)=0$ et comme $m$ appartient à $B(\widetilde{\ell}'_0,\widetilde{\ell})\subset 
F\cap(B(\widetilde{\ell}'_{-1},\widetilde{\ell}_0)$, il
est entrelacé avec $\widetilde{\alpha}_\gamma$ par hypothèse.
Donc l'image de $\widetilde{\ell}'_0$ dans $T$ appartient
à l'axe de 
translation de $\gamma$ dans $T$ et par $\gamma$-invariance, l'image de $\widetilde{\ell'}$ aussi. Donc si $i_{[q]}(\mu,\alpha_\gamma)>0$, la distance de 
translation $\ell_T(\gamma)$ est égale à $i_{[q]}(\mu,\alpha_\gamma)$ et l'axe de translation de $\gamma$
est l'image de l'ensemble des feuilles de $\widetilde{\Lambda}$ qui sont entrelacées avec un axe de 
translation de $\gamma$ sur $\revet$.

{\bf Cas $(2)$.~} Supposons que $i_{[q]}(\mu,\alpha_\gamma)=0$, c'est-à-dire que $\widetilde{\alpha}_\gamma$ n'est entrelacée avec aucune 
feuille de $\widetilde{\Lambda}$. Si $\widetilde{\alpha}_\gamma$ a le même couple de points à l'infini qu'une feuille $\widetilde{\ell}$ de 
$\widetilde{\Lambda}$,
alors $\gamma\widetilde{\ell}=\widetilde{\ell}$ et $\widetilde{\ell}^T$ est un point fixe de $\gamma$ dans $T$.
Sinon, on note
$(S,N)=(\widetilde{\alpha}_\gamma(-\infty),\widetilde{\alpha}_\gamma(+\infty))$. D'après le lemme \ref{quedescomposnatesminimales},
aucune feuille de $\Lambda$ n'est positivement périodique sans être périodique, donc d'après \cite[Lem.~4.13~et~4.14]{Morzy1}, les points $N$ et $S$
ne sont un point à l'infini d'aucune feuille de $\widetilde{\Lambda}$.
On rappelle que l'ordre cyclique total $o$ sur $\partial_\infty\widetilde{\Sigma}$ 
définit une relation d'ordre total $\leqslant$ sur $\partial_\infty\widetilde{\Sigma}$ défini par $S<\eta$ pour tout 
$\eta\in\partial_\infty\widetilde{\Sigma}-\{S\}$ et $\eta_1\leqslant\eta_2$ si et seulement si
$o(\eta_1,\eta_2,S)\in\{0,1\}$ pour tout $\eta_1,\eta_2\in\partial_\infty\widetilde{\Sigma}-\{S\}$ 
(voir \cite[Rem.~2.9]{Morzy1} pour la définition de $o$ et \cite[Déf.~2.23]{Wolf11} pour la définition de $\leqslant$).
Soient $\widetilde{\ell}$ une feuille de $\widetilde{\Lambda}$ et $(a,b)=(\widetilde{\ell}(-\infty),\widetilde{\ell}(+\infty))$ son couple de points à l'infini.
 Puisque les feuilles de
$\widetilde{\Lambda}$ sont deux-à-deux non entrelacées, si 
$(a',b')\in\dddp$ est le couple de points à l'infini d'une feuille $\widetilde{\ell}'$, telle que $S\leqslant a'\leqslant a$ et $a\leqslant b'\leqslant N$, 
alors $b\leqslant b'\leqslant N$.
Donc,
  par compacité de $\partial_\infty\widetilde{\Sigma}$, quitte à remplacer
$\widetilde{\ell}$, on peut supposer qu'il n'existe pas de telle feuille $\widetilde{\ell}'$ telles que $S< a'< a$ et $b< b'< N$. De même, si 
$\widetilde{\ell}$ est contenue dans une bande plate, on peut toujours supposer qu'il n'existe pas de feuille contenue dans cette bande plate qui 
soit contenue dans l'intersection des adhérences des composantes connexes de $\widetilde{\Sigma}-\widetilde{\ell}(\RR)$ contenant $\widetilde{\alpha}_\gamma$
et de  $\widetilde{\Sigma}-\widetilde{\alpha}_\gamma(\RR)$ 
contenant $\widetilde{\ell}$, hormis $\widetilde{\ell}$. 
Or, l'action de $\gamma$ sur $\partial_\infty\widetilde{\Sigma}$ est une action de type Nord-Sud, dont les points fixes sont $N$ et $S$,
représentée par les flèches sur le dessin.

\medskip
\noindent
\begin{minipage}{10 cm}
Donc $a<\gamma a$, $b<\gamma b$ et puisque $\widetilde{\ell}$ et $\gamma\widetilde{\ell}$ ne sont pas entrelacées, on a aussi 
$b\leqslant \gamma a\leqslant \gamma b$. Supposons qu'il existe une feuille $\widetilde{\ell}'\in B(\widetilde{\ell},\gamma\widetilde{\ell})
-\{\widetilde{\ell},\gamma\widetilde{\ell}\}$,
de couple de points à l'infini $(a',b')$. Comme $\widetilde{\ell}'$ n'est pas entrelacée avec $\widetilde{\alpha}_\gamma$ et par hypothèse sur 
$\widetilde{\ell}$, quitte à remplacer $\widetilde{\ell}'$ par son inverse, on  a $b\leqslant a'\leqslant \gamma a$ et 
$\gamma b\leqslant b'< N$, et $(a',b')\not=(\gamma a,\gamma b)$. Mais alors $S<\gamma^{-1}a'\leqslant a$ et $b\leqslant\gamma^{-1}b'< N$, avec
$(\gamma^{-1}a',\gamma^{-1}b')\not=(a,b)$, ce qui contredit l'hypothèse faite sur $\widetilde{\ell}$. Donc
$B(\widetilde{\ell},\gamma\widetilde{\ell})=\{\widetilde{\ell},\gamma\widetilde{\ell}\}$ et $\widetilde{d}(\widetilde{\ell},\gamma\widetilde{\ell})=0$.
Donc $\widetilde{\ell}^T$ est un point fixe de $\gamma$ dans $T$.\cqfd
\end{minipage}
\begin{minipage}{4.9 cm}
\begin{picture}(0,0)%
\includegraphics{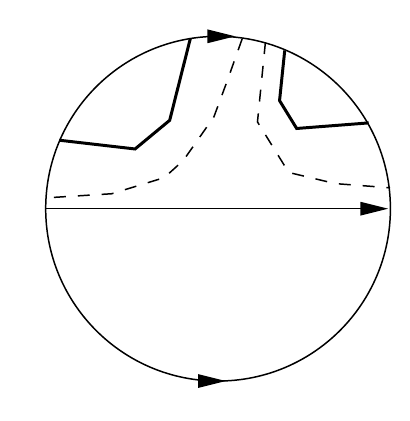}%
\end{picture}%
\setlength{\unitlength}{3771sp}%
\begingroup\makeatletter\ifx\SetFigFont\undefined%
\gdef\SetFigFont#1#2#3#4#5{%
  \reset@font\fontsize{#1}{#2pt}%
  \fontfamily{#3}\fontseries{#4}\fontshape{#5}%
  \selectfont}%
\fi\endgroup%
\begin{picture}(2050,2107)(-1094,-224)
\put(-455,579){\makebox(0,0)[lb]{\smash{{\SetFigFont{11}{13.2}{\familydefault}{\mddefault}{\updefault}{\color[rgb]{0,0,0}$\widetilde{\alpha}_\gamma$}%
}}}}
\put(-984,1171){\makebox(0,0)[lb]{\smash{{\SetFigFont{11}{13.2}{\familydefault}{\mddefault}{\updefault}{\color[rgb]{0,0,0}$a$}%
}}}}
\put(941,766){\makebox(0,0)[lb]{\smash{{\SetFigFont{11}{13.2}{\familydefault}{\mddefault}{\updefault}{\color[rgb]{0,0,0}$N$}%
}}}}
\put(-233,1736){\makebox(0,0)[lb]{\smash{{\SetFigFont{11}{13.2}{\familydefault}{\mddefault}{\updefault}{\color[rgb]{0,0,0}$b$}%
}}}}
\put(784,1234){\makebox(0,0)[lb]{\smash{{\SetFigFont{11}{13.2}{\familydefault}{\mddefault}{\updefault}{\color[rgb]{0,0,0}$\gamma_b$}%
}}}}
\put(317,1695){\makebox(0,0)[lb]{\smash{{\SetFigFont{11}{13.2}{\familydefault}{\mddefault}{\updefault}{\color[rgb]{0,0,0}$\gamma a$}%
}}}}
\put(-107,156){\makebox(0,0)[lb]{\smash{{\SetFigFont{11}{13.2}{\familydefault}{\mddefault}{\updefault}{\color[rgb]{0,0,0}$\widetilde{\Sigma}$}%
}}}}
\put(-1079,794){\makebox(0,0)[lb]{\smash{{\SetFigFont{11}{13.2}{\familydefault}{\mddefault}{\updefault}{\color[rgb]{0,0,0}$S$}%
}}}}
\put(721,1424){\makebox(0,0)[lb]{\smash{{\SetFigFont{11}{13.2}{\familydefault}{\mddefault}{\updefault}{\color[rgb]{0,0,0}$\gamma\widetilde{\ell}$}%
}}}}
\put(406,-151){\makebox(0,0)[lb]{\smash{{\SetFigFont{11}{13.2}{\familydefault}{\mddefault}{\updefault}{\color[rgb]{0,0,0}$\partial_\infty\widetilde{\Sigma}$}%
}}}}
\put(-809,1469){\makebox(0,0)[lb]{\smash{{\SetFigFont{11}{13.2}{\familydefault}{\mddefault}{\updefault}{\color[rgb]{0,0,0}$\widetilde{\ell}$}%
}}}}
\put(-179,929){\makebox(0,0)[lb]{\smash{{\SetFigFont{11}{13.2}{\familydefault}{\mddefault}{\updefault}{\color[rgb]{0,0,0}$\gamma^{-1}\widetilde{\ell}'$}%
}}}}
\put(496,1019){\makebox(0,0)[lb]{\smash{{\SetFigFont{11}{13.2}{\familydefault}{\mddefault}{\updefault}{\color[rgb]{0,0,0}$\widetilde{\ell}'$}%
}}}}
\end{picture}%

\end{minipage}

\section{Liens entre l'arbre dual à une lamination plate mesurée et l'arbre dual à la lamination hyperbolique mesurée correspondante.}

Dans cette partie, on reprend les notations et conventions de la partie \ref{liensplatshyperboliques}. Nous commençons par rappeler la définition 
de l'arbre dual à une lamination hyperbolique mesurée (voir par exemple \cite[§1]{MorSha91}), avec une présentation nouvelle qui nous permettra de définir une
isométrie $\grperevet$-équivariante entre l'arbre dual à une lamination plate mesurée et l'arbre dual à la lamination hyperbolique mesurée correspondante.
Soit $(\Lambda_m,\mu_m)$ une lamination hyperbolique
mesurée de $\srfcem$ et $(\Lmr,\mutilde_{{m}})$ son image réciproque dans $\widetilde{\Sigma}$. Alors $(\Lmr,\mutilde_{{m}})$ est $\grperevet$-invariante et elle définit une mesure 
$\nu_{\mutilde_m}\in\M_{\grperevet}([\Gmr])$ (voir \cite[Prop.~17~p.~154]{Bonahon88}). Si $\widetilde{\lambda}$ est un atome de $\nu_{\mutilde_m}$, on remplace 
$\widetilde{\lambda}$ par une bande plate de largeur $\nu_{\mutilde_m}(\widetilde{\lambda})$ feuilletée par des droites parallèles aux bords. En procédant ainsi pour chacun des
atomes de $\nu_{\mutilde_m}$, on obtient une lamination géodésique $\widetilde{\Lambda}'$ sur une surface munie d'une métrique $(\widetilde{\Sigma}',
d')$, qui est $\CAT(0)$ et localement
$\CAT(-1)$ dans le complémentaire des bandes plates, et l'action par isométries de $\grperevet$ sur $\revetm$ se prolonge de manière unique en une action
par isométries sur $(\widetilde{\Sigma}',d')$. 
     Remarquons que l'espace $(\widetilde{\Sigma}',d')$ est un espace métrique $\CAT(0)$ enrubanné, et on est dans le cadre de la généralisation des
     laminations géodésiques
proposée dans \cite{Morzy1}.

De plus, pour chaque bande plate $\BP(\widetilde{\lambda})$ associée à un atome $\widetilde{\lambda}$ de $\nu_{\mutilde_m}$, si
$\alpha:[0,\nu_{\mutilde_m}(\widetilde{\lambda})]\to \BP(\widetilde{\lambda})$ est un segment qui 
relie orthogonalement les deux bords de $\BP(\widetilde{\lambda})$, alors l'application $r_\alpha$ de l'ensemble des feuilles de $\widetilde{\Lambda}'$ contenues dans $\BP(\widetilde{\lambda})$
dans $\Image(\alpha)$ définie par $r(g)=g(\RR)\cap\Image(\alpha)$ est un homéomorphisme. Donc on peut munir l'ensemble des feuilles contenues dans $\BP(\widetilde{\lambda})$
de la mesure $\nu_{\widetilde{\lambda}}=(r_{\alpha}^{-1})_*dx_ \alpha$ où 
$dx_\alpha$ est la mesure proportionnelle à la mesure de Lebesgue sur $\Image(\alpha)$, de masse $\nu_{\widetilde{\mu}_m}(\widetilde{\lambda})$. 
La mesure $\nu'$ égale à $\nu_{\mutilde_m}$ en dehors des atomes et égale à $\nu_{\widetilde{\lambda}}$ sur l'ensemble des feuilles
feuilletant la bande plate associée à $\widetilde{\lambda}$, pour chaque atome $\widetilde{\lambda}$ de $\nu_{\mutilde_m}$, 
est une mesure de Radon sur $[\G_{d'}]$ qui est 
$\grperevet$-invariante, sans atome et de support égal à $\widetilde{\Lambda}'$. 

Si $\widetilde{\lambda}_0$ et $\widetilde{\lambda}_1$ sont des feuilles de $\widetilde{\Lambda}'$, 
et si $c$ est un segment géodésique joignant leurs images, on note $B_{\widetilde{\Lambda}'}(\widetilde{\lambda}_0,\widetilde{\lambda}_1)$ (ou plus rapidement 
$B(\widetilde{\lambda}_0,\widetilde{\lambda}_1)$) l'ensemble compact des feuilles de 
$\widetilde{\Lambda}'$ qui sont contenues dans l'intersection des adhérences des composantes connexes des complémentaires de 
$\widetilde{\lambda}_0(\RR)$ et $\widetilde{\lambda}_1(\RR)$ qui contiennent respectivement $\widetilde{\lambda}_1$ et $\widetilde{\lambda}_0$, 
et qui intersectent $c$ non trivialement. Alors de même que dans la partie \ref{arbredual}, 
l'ensemble $B_{\widetilde{\Lambda}'}(\widetilde{\lambda}_0,\widetilde{\lambda}_1)$ ne dépend pas du choix de $c$, il est compact et muni d'un ordre total, que l'on note
toujours $\preceq$, défini par $\widetilde{\lambda}\preceq\widetilde{\lambda}'$ si et seulement si $\widetilde{\lambda}\in B(\widetilde{\lambda}_0,
\widetilde{\lambda}')$. On définit aussi une pseudo-distance
$\widetilde{d}_{\widetilde{\Lambda}'}$ sur $\widetilde{\Lambda}'$ par
$\widetilde{d}_{\widetilde{\Lambda}'}(\widetilde{\lambda}_0,\widetilde{\lambda}_1)=\nu_{\mutilde_m}(B_{\widetilde{\Lambda}'}(\widetilde{\lambda}_0,\widetilde{\lambda}_1))$
pour toutes les feuilles $\widetilde{\lambda}_0$, $\widetilde{\lambda}_1$ de $\widetilde{\Lambda}'$, et 
le quotient de $(\widetilde{\Lambda}',d_{\widetilde{\Lambda}'})$ par la relation d'équivalence 
$\widetilde{\lambda}_0\sim\widetilde{\lambda}_1$ si et seulement si 
$d_{\widetilde{\Lambda}'}(\widetilde{\lambda}_0,\widetilde{\lambda}_1)=0$ (ou de manière équivalente, 
$B(\widetilde{\lambda}_0,\widetilde{\lambda}_1)=\{\widetilde{\lambda}_0,\widetilde{\lambda}_1\}$) est un arbre réel $(T,d_T)$ que l'on appelle {\it l'arbre dual} à $(\Lmr,\mutilde_m)$. 
Pour tout 
$\gamma\in\grperevet$,
$\gamma_*\nu'=\nu'$ et 
$\gamma B(\widetilde{\lambda}_0,\widetilde{\lambda}_1)=B(\gamma\widetilde{\lambda}_0,\gamma\widetilde{\lambda}_1)$ pour toutes
les feuilles $\widetilde{\lambda}_0$ et $\widetilde{\lambda}_1$ de $\widetilde{\Lambda}'$, donc l'action de $\grperevet$ sur $\Lmr$ passe au quotient et définit
une action par isométries sur $(T,d_T)$.
On se convainct aisément qu'il existe une isométrie $\grperevet$-équivariante entre l'arbre dual ainsi construit, et celui construit par exemple dans 
\cite[§1]{MorSha91}, en identifiant les feuilles de $\widetilde{\Lambda}'$ et les composantes connexes du complémentaire du support de $\widetilde{\Lambda}'$
qu'elles bordent.

\medskip

Soient $(\Lqr,\mutilde_{[q]})$ une lamination plate mesurée de $\revet$, $\nu_{\mutilde_{[{q}]}}$ la mesure qu'elle définit sur $[\G_{[\widetilde{q}]}]$
et $\nu_{\mutilde_m}$ son image par $\varphi_*$ (voir le lemme \ref{surjectionpropre}). On note
$(\Lmr,\mutilde_{m})$ la lamination hyperbolique mesurée définie par $\nu_{\mutilde_m}$, et $\widetilde{\Lambda}'$ et $\nu'$ les lamination géodésique 
sur $(\widetilde{\Sigma}',d')$ et mesure de Radon sur $[\G_{d'}]$ associées à $(\Lmr,\mutilde_m)$ par la construction ci-dessus.
On suppose (quitte à procéder comme dans
le dernier paragraphe de la partie \ref{arbredual}) que $\nu_{\mutilde_{[{q}]}}$ n'a pas d'atome.

Si $\widetilde{\lambda}$ est un atome de $\nu_{\mutilde_m}$, et si $F_{\widetilde{\lambda}}$ est l'ensemble des feuilles de $\Lqr$ auxquelles correspond 
$\widetilde{\lambda}$ 
(voir \cite[§4.2]{Morzy1}), 
il existe des bandes plates maximales de $\revet$ et de $(\widetilde{\Sigma}',d')$ qui contiennent respectivement $F_{\widetilde{\lambda}}$ et 
l'ensemble $F'_{\widetilde{\lambda}}$ des feuilles de
$\widetilde{\Lambda}'$ qui correspondent à $\widetilde{\lambda}$ dans la construction ci-dessus. On note $\widetilde{\ell}_0$ et $\widetilde{\ell}_1$ (resp.
$\widetilde{\lambda}_0$ et $\widetilde{\lambda}_1$) les feuilles extrémales de $F_{\widetilde{\lambda}}$ (resp. $F'_{\widetilde{\lambda}}$), c'est-à-dire telles 
que $F_{\widetilde{\lambda}}=B_{\Lmr}(\widetilde{\ell}_0,\widetilde{\ell}_1)$
et $F_{\widetilde{\lambda}}'=B_{\widetilde{\Lambda}'}(\widetilde{\lambda}_0,\widetilde{\lambda}_1)$.
Alors il existe une unique application $\phi_{\widetilde{\lambda}}
:F_{\widetilde{\lambda}}\to F_{\widetilde{\lambda}}'$ telle que pour toute 
$\widetilde{\ell}\in F_{\widetilde{\lambda}}$, on ait $\nu_{\mutilde_m}(B_{\widetilde{\Lambda}'}(\phi_{\widetilde{\lambda}}(\widetilde{\ell}_0),
\phi_{\widetilde{\lambda}}(\widetilde{\ell}))=
\nu_{\mutilde_{[q]}}(B_{\Lqr}(\widetilde{\ell}_0,\widetilde{\ell}))$.
On note alors $\phi$ l'application de $\Lqr$ dans $\widetilde{\Lambda}'$
égale à $\varphi_{|\Lqr}$ en dehors des images réciproques des atomes de $\nu_{\mutilde_m}$ (voir \cite[§4.2]{Morzy1}) et à $\phi_{\widetilde{\lambda}}$ 
sur les ensembles 
$F_{\widetilde{\lambda}}$ où $\widetilde{\lambda}$ est un atome de $\nu_{\mutilde_m}$. Alors par construction, pour toutes les 
feuilles $\widetilde{\ell}_0$ et  $\widetilde{\ell}_1$ de $\Lqr$, on a $B_{\widetilde{\Lambda}'}(\phi(\widetilde{\ell}_0),\phi(\widetilde{\ell}_1))=
\phi(B_{\Lqr}(\widetilde{\ell}_0,\widetilde{\ell}_1))$. Donc l'application $\phi$ passe au quotient et définit 
une application 
$\phi_T:(T_{[q]},d_{T_{[q]}})\to(T_{m},d_{T_{m}})$, où 
$(T_{[q]},d_{T_{[q]}})$ et $(T_{m},d_{T_{m}})$ sont respectivement les arbres duaux à
$(\Lqr,\mutilde_{[q]})$ et à $(\Lmr,\mutilde_m)$.

\blemm \label{isometrieequivariante}
L'application $\phi_T:(T_{[q]},d_{T_{[q]}})\to(T_{m},d_{T_{m}})$ est une isométrie $\grperevet$-équivariante.
\elemm

\dem Si $\widetilde{\ell}_0$ et $\widetilde{\ell}_1$ sont des feuilles de $\Lqr$, on a
$\phi(B_{\Lqr}(\widetilde{\ell}_0,\widetilde{\ell}_1))=B_{\widetilde{\Lambda}'}(\phi(\widetilde{\ell}_0),\phi(\widetilde{\ell}_1))$
et comme 
$\nu_{\mutilde_m}=\varphi_*\nu_{\mutilde_{[q]}}$, l'application $\phi_T$ est isométrique. De plus, $\phi$ est surjective donc, par passage au quotient,
$\phi_T$ aussi et c'est une isométrie. Enfin, comme $\varphi$, l'application
$\phi$ est $\grperevet$-équivariante et par passage au quotient, $\phi_T$ aussi.\cqfd

%
%

\bibliographystyle{alphanum}
\bibliography{biblio}{}
Département de mathématique, UMR 8628 CNRS, Université Paris-Sud, Bât. 430, F-91405 Orsay Cedex, FRANCE. 
Bureau : 16.

{\it thomas.morzadec@math.u-psud.fr}
\end{document}